\documentclass[10pt]{amsart}

\usepackage{tikz}
\usetikzlibrary{positioning}
\usetikzlibrary{arrows}
\usepackage{amssymb,amsmath,latexsym,graphicx}
\usepackage{charter,eucal}
\usepackage{amssymb}
\usepackage{amscd}
\usepackage[all,cmtip]{xy}
\usepackage{rotating,multicol,soul}
\usetikzlibrary{decorations.markings}

\setul{}{1.1pt}

\newtheorem{theorem}{Theorem }[section]

\newtheorem{lemma}[theorem]{Lemma}
\newtheorem{observation}[theorem]{Observation}

\newtheorem{remark}[theorem]{Remark}
\newtheorem{corollary}[theorem]{Corollary}
\newtheorem{proposition}[theorem]{Proposition}
\newtheorem{principle}[theorem]{\textsc{Principle}}

\newcommand{\bt}{\begin{theorem}}
\newcommand{\et}{\end{theorem}}
\newcommand{\bmt}{\begin{maintheorem}}
\newcommand{\emt}{\end{maintheorem}}
\newcommand{\bc}{\begin{corollary}}
\newcommand{\bl}{\begin{lemma}}
\newcommand{\ec}{\end{corollary}}
\newcommand{\el}{\end{lemma}}
\newcommand{\bo}{\begin{observation}}
\newcommand{\eo}{\end{observation}}
\newcommand{\bp}{\begin{proposition}}
\newcommand{\ep}{\end{proposition}}
\newcommand{\br}{\begin{remark}}
\newcommand{\er}{\end{remark}}
\newcommand{\bpr}{\begin{principle}}
\newcommand{\epr}{\end{principle}}

\newcommand{\hooklongrightarrow}{\lhook\joinrel\longrightarrow}

\newcommand{\comm}[1]{{\color{red}{\bf #1}}}

\def\Aut{\mathrm{Aut}}
\def\I{\mathop{\mathrm{I}}}

\def\PG{\mathbf{PG}}

\def\I{\mathbf{I}}

\def\eop{\hspace*{\fill}$\blacksquare$}

\def\id{\mathrm{id}}

\newcommand{\F}{\mathbb{F}}

\newcommand{\mC}{\mathcal{C}}
\newcommand{\mG}{\mathcal{G}}

\newcommand{\mQ}{\mathcal{Q}}

\newcommand{\mA}{\mathcal{A}}
\newcommand{\mP}{\mathcal{P}}
\newcommand{\mL}{\mathcal{L}}
\newcommand{\mE}{\mathcal{E}}

\newcommand{\mW}{\mathcal{W}}

\newcommand{\mO}{\mathcal{O}}

\newcommand{\mS}{\mathcal{S}}

\newcommand{\mH}{\mathcal{H}}

\newcommand{\mK}{\mathcal{K}}
\newcommand{\hT}{\mathbf{T}}

\newcommand{\Hom}{\texttt{hom}}
\newcommand{\cor}{\texttt{cor}}

\usetikzlibrary{matrix}
\usetikzlibrary{arrows}
\usepackage{pgf}
\usepackage{lscape}
\usepackage{listings}

\usepackage{enumitem}

\title[Covers of quadrangles, 2]{Covers of generalized quadrangles, 2. Kantor-Knuth covers and embedded ovoids}



\author{Joseph A. Thas and Koen Thas}

\thanks{}

\address{Ghent University, Department of Mathematics, Krijgslaan 281, S22 and S25, B-9000 Ghent, Belgium}

\email{thas.joseph@gmail.com; koen.thas@gmail.com}

\date{}

\begin{document}

\maketitle

\begin{abstract}
In this paper, which is a sequel to \cite{part1}, we proceed with our study of covers and decomposition laws for geometries related to generalized quadrangles. In particular, we obtain a higher decomposition law for all Kantor-Knuth generalized quadrangles which generalizes one of the main results in \cite{part1}. In a second part of the paper, we study the set of all Kantor-Knuth ovoids (with given parameter)  
in a fixed finite parabolic quadrangle, and relate this set to embeddings of parabolic quadrangles into Kantor-Knuth quadrangles. This point of view gives rise to an answer of a question posed in \cite{JATSEP}.  
\end{abstract}

\setcounter{tocdepth}{1}
\bigskip
{\footnotesize
\tableofcontents
}

\medskip
\section{Introduction}

\subsection{Decomposition laws}

In the paper \cite{CS}, Cardinali and Sastry study various aspects of geometries related to ovoids of the classical
finite symplectic generalized quadrangle $\mW(2^n)$. In particular, starting with the natural embedding $\mW(2^n) \hookrightarrow \mQ(5,2^n)$, 
they consider the geometry $\mA$ defined by $\mQ(5,2^n) \setminus \mW(2^n)$, and the geometry $\mE$ which consists of elliptic 
ovoids in $\mW(2^n)$ as points, and rosettes of such ovoids as lines (details related to these notions can be found in the next two sections). There is a natural
projection 
\begin{equation}
\pi: \mA \mapsto \mE.
\end{equation} 

Cardinali and Sastry ask if any given $2$-cover $\gamma: \mA \mapsto \mE$ factors through $\pi$ and an automorphism $\widetilde{\alpha}$
of $\mA$ (or $\mQ(5,2^n)$):

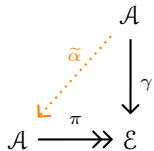
\bigskip
\begin{center}
\begin{tikzpicture}[>=angle 90,scale=2.2,text height=2.0ex, text depth=0.45ex]
\node (a0) at (0,3) {$\mA$};
\node (a1) [right=of a0] {$\mE$};

\node (b0) [above=of a1] {$\mA$};

\draw[->>,font=\scriptsize,thick]

(a0) edge node[above] {$\pi$} (a1);

\draw[->,font=\scriptsize,thick]

(b0) edge node[right] {$\gamma$} (a1);

\draw[->,font=\scriptsize,blue,thick,dotted,orange]
(b0) edge node[above] {$\widetilde{\alpha}$} (a0);

\end{tikzpicture}
\end{center}


That question was the starting point of \cite{part1}, and in a more general setting |  without finiteness
restrictions, and without assuming that the quadrangles in question are classical | a complete answer was provided in loc. cit. 

More precisely, we have:

\bt[\cite{part1}]
\label{MT1}
Let $\mS$ be a thick generalized quadrangle, and let $\mS' \hookrightarrow \mS$ be a thick full subquadrangle which is a geometrical hyperplane of $\mS$. 
Define $\mA$ and $\mE$ as above, and let $\pi: \mA \mapsto \mE$ be the natural projection.
Then every cover $\gamma: \mA \mapsto \mE$ factorizes as $\gamma = \pi \circ \widetilde{\alpha}$, with $\widetilde{\alpha}$ an automorphism of $\mA$, if and only if every automorphism of $\mE$ is induced by an automorphism of $\mS$. 
\et

In \cite{part1}, we called the factorization  $\gamma = \pi \circ \widetilde{\alpha}$ ``higher decomposition,'' and in \cite{part1}, other factorizations were studied (cf. section \ref{decomps} for more details). Plugging in the data $\Big(\mS, \mS' \Big) = \Big(\mQ(5,q),\mW(q) \Big)$  in Theorem \ref{MT1} readily leads to the answer of the question 
of Cardinali and Sastry for all prime powers $q$.  

The first main result of the present paper yields higher decomposition for the class of finite Kantor-Knuth quadrangles (with parabolic subquadrangles). For $q$ an odd prime power, this generalizes our answer to the question of Cardinali and Sastry, since $\mQ(5,q)$ is just a classical Kantor-Knuth generalized quadrangle. 

In particular, we will obtain the following result in section \ref{s2} (where $\mA$ and $\mE$ are constructed from a nonclassical Kantor-Knuth quadrangle with 
a full parabolic subquadrangle):

\begin{theorem}
Let $\gamma: \mathcal{A} \longmapsto \mathcal{E}$ be a 2-cover of $\mathcal{E}$. Then $\gamma$ induces an automorphism $\overline \alpha$ of $\mathcal{S^\prime}$ and there exist exactly two automorphisms $\widetilde \alpha$ of $\mathcal{S}$ which extend the automorphism $\overline \alpha$ of $\mathcal{S^\prime}$.
\end{theorem}

\subsection{A garden of Kantor-Knuth ovoids}

The starting point of the second part of this paper is the following quote (taken from \cite{JATSEP}, pp. 245):

\begin{quote}
``Note that the number $q^2(q^2 - 1)(q^4 - 1)/4$ of ovoids produced for $q$ a power of $3$ is much larger than the number of points 
$z \in \mS \setminus \mS_{\gamma}$ available to be used in the original construction.''
\end{quote}

Indeed, although Kantor-Knuth ovoids are initially defined as subtended ovoids $\mO_x$ in parabolic quadrangles $\mQ(4,q)$ which are 
full subquadrangles in Kantor-Knuth quadrangles, a calculation shows that the $\Aut(\mQ(4,q))$-orbit of $\mO_x$ contains much more elements than the number 
of ovoids obtained as subtended ovoids. In section \ref{emb}, we try to understand this difference, by observing that one can partition the orbit in certain equivalence classes, each class of which corresponds to an embedding of $\mQ(4,q)$ in a Kantor-Knuth quadrangle (of the same type as the Kantor-Knuth ovoid). 

In order to facilitate this view, we study a category of embeddings of parabolic quadrangles in Kantor-Knuth quadrangles in much detail, and meet a number of interesting combinatorial problems related to covers on the wayside.

\subsection{Organization} 

In section \ref{basic}, we sum up some basic notions concerning the geometry of generalized quadrangles (and related geometries). In section \ref{decomps}, we 
recall higher and lower decomposition from \cite{part1}. Sections \ref{recap} and \ref{tranov} provide a short overview of the basic  properties 
of  Kantor-Knuth quadrangles, Kantor-Knuth ovoids and translation ovoids which we will need in our paper. Then, in section \ref{s2}, higher decomposition is obtained 
for all Kantor-Knuth quadrangles. In section \ref{emb}, we study embeddings of parabolic quadrangles in Kantor-Knuth quadrangles and the various related combinatorial problems which arise. Finally, in the last section, we calculate which fraction of the size of the automorphism group of a Kantor-Knuth ovoid is induced 
by the automorphism group of a Kantor-Knuth quadrangle.

Our approach in this paper is mainly synthetic. 

\section{Basic geometric notions}
\label{basic}

We summarize some basic notions that will be used throughout. Standard references are \cite{HT2,PT2,TGQ}.

\subsection{Quadrangles}

A {\em generalized quadrangle} (GQ) of {\em order $(s,t)$}, where $s$ and $t$ are cardinal numbers,
is a point-line incidence geometry with the following properties:
\begin{itemize}
\item[(a)]
Each point is incident with $t + 1$ lines and each line is incident with $s + 1$ points.
\item[(b)]
If $(x,L)$ is a non-incident point-line pair, there is precisely one point-line pair $(y,M)$ such that
$x \I M \I y \I L$ (here, ``$\I$'' denotes the incidence relation).
\item[(c)]
Two distinct points are incident with at most one line.
\end{itemize}

If both $s$ and $t$ are at least $2$, we say that the GQ is {\em thick}; otherwise it is called {\em thin}.

We use the usual notation $x \sim y$ to indicate that the points $x$ and $y$ are collinear, and any point is collinear with itself. For any point $x$, $x^{\perp} := \{ z \vert z \sim x \}$, and for any subset $Y$ of the point set, $Y^{\perp} := \cap_{y \in Y}y^{\perp}$; we make no distinction between $y^{\perp}$ and $\{y\}^{\perp}$.
In particular, we denote ${(Y^{\perp})}^{\perp}$ by $Y^{\perp\perp}$.

Finally, if $u$ and $v$ are distinct points, then $\mathrm{cl}(u,v)$ is the point set $\{ w \vert w^{\perp} \cap \{u,v\}^{\perp\perp} \ne \emptyset \}$.   

\subsection{Subquadrangles}

If $\mS = (\mP,\mL,\I)$ is a generalized quadrangle (where $\mP$ is the point set, $\mL$ the line set, and $\I$ is incidence), then a 
{\em subgeometry} is a triple $\Gamma = (\mP',\mL',\I')$ such that 
$\mP' \subseteq \mP$, $\mL' \subseteq \mL$, and $\I'$ is the induced incidence relation. 

A {\em subquadrangle} (subGQ) is a subgeometry which is a generalized quadrangle.
A subquadrangle $\mS'$ of a generalized quadrangle $\mS$ is {\em full} if for all its lines $U$, we have the property that 
any point of $\mS$ incident with $U$ is also a point of $\mS'$. Dually, we speak of ``ideal subquadrangles.''

\subsection{Ovoids}

If $\mS$ is a generalized quadrangle, an {\em ovoid} of $\mS$ is a set of points $\mO$ in $\mS$ such that each line 
contains precisely one point of $\mO$.

\subsection{Semi partial geometries}

A (finite) {\em semi partial geometry} is an incidence structure $\mathcal{I} = (\mathcal{P}, \mathcal{B}, {\rm I})$, in which $\mathcal{P}$ and $\mathcal{B}$ are disjoint, non-empty sets of points and lines, for which ${\rm I}$ is a symmetric point-line incidence relation satisfying the following axioms:
\begin{itemize}
\item [{\rm (i)}] each point is incident with $1 + t^\ast$ lines, with $t^\ast \geq 1$, and two distinct points are incident with at most one line; 
\item [{\rm (ii)}] each line is incident with $1 + s^\ast$ points, with $s^\ast \geq 1$, and two distinct lines are incident with at most one point;
\item [{\rm (iii)}] if a point $x$ and a line $L$ are not incident, then there are 0 or $\alpha^\ast$ points, where $\alpha^\ast \geq 1$, which are collinear with $x$ and incident with $L$;
\item [{\rm (iv)}] if two points are not collinear, then there are $\mu^\ast$ points, where $\mu^\ast > 0$, collinear with both.
\end{itemize}

The integers $s^\ast, t^\ast, \alpha^\ast, \mu^\ast$ are the {\em parameters} of the semi partial geometry. The semi partial geometries with $\alpha^\ast = 1$ are the {\em partial quadrangles}. A semi partial geometry is a {\em partial geometry} if and only if the 0 in Axiom (iii) does not occur; this is equivalent to the condition $\mu^\ast = (t^\ast + 1)\alpha^\ast$.

A {\em generalized quadrangle} is a partial geometry with $\alpha^\ast = 1$, so $\mu^\ast = t^\ast + 1$.

For more on semi partial geometries we refer to \cite{HT2}. 

\subsection{Morphisms}

If $\Gamma = (\mP,\mL,\I)$ and $\Gamma' = (\mP',\mL',\I')$ are point-line incidence geometries, then a {\em morphism}
$\gamma: \Gamma \mapsto \Gamma'$ is a map from $\mP \cup \mL$ to $\mP' \cup \mL'$ which sends points to points, lines to lines, and preserves incidence. If $\gamma$ is bijective and the inverse map also preserves incidence, we say $\gamma$ is an {\em isomorphism}. If $\Gamma = \Gamma'$, an isomorphism is also called an {\em automorphism}, and the set of automorphisms of $\Gamma$ naturally forms a group under composition of maps, which is denoted by $\Aut(\Gamma)$.

\section{Higher and lower decomposition}
\label{decomps}

In this section, we describe higher and lower decomposition as introduced in \cite{part1}.

\subsection {$\mathcal{A}$ and $\mathcal{E}$}\label{intro1} 

Let $\mathcal{S}$ be a thick generalized quadrangle, and let $\mathcal{S^\prime}$ be a thick full subquadrangle of $\mathcal{S}$. Let $x$ be a point of $\mathcal{S} \setminus 
\mathcal{S^\prime}$. Then $x^\perp \cap \mathcal{S^\prime} = \mO_x$ is easily seen to be an ovoid of $\mathcal{S^\prime}$, which we call a {\em{subtended ovoid (by $x$)}}. Let $L$ be any line of $\mathcal{S} \setminus \mathcal{S^\prime}$ which meets $\mathcal{S^\prime}$ in a point $l$. Then each point of $L \setminus \lbrace l \rbrace$ subtends an ovoid in $\mathcal{S^\prime}$, and the set of all these ovoids is called the {\em{rosette}} $R_L$ of ovoids determined by $L$. If $\mathcal{S^\prime}$ is also a geometrical hyperplane, there is a natural bijection between the lines incident with $x$ and the points of $x^\perp \cap \mathcal{S^\prime}$. 

Keeping the latter hypothesis, we define the geometry $\mathcal{E}$ to have as points the subtended ovoids in $\mathcal{S^\prime}$, and as lines the rosettes $R_L$. Incidence is symmetrized containment. We also define a geometry $\mathcal{A}$, which is just the affine generalized quadrangle which arises when taking away the geometric hyperplane $\mathcal{S^\prime}$ of $\mathcal{S}$. 

Note that there is a natural projection \\
 \begin{equation}
 \pi : \mathcal{A} \longmapsto \mathcal{E} , \left \{ \begin{array}{rl} x \longmapsto \mO_x & \mbox{for all points $x$} \\ L \longmapsto R_L & \mbox{for all lines $L$.} \end{array} \right . 
\end{equation}

By the mere definition of $\mathcal{E}$, $\pi$ is surjective on both points and lines.

\subsection{Covers}\label{intro2} 

Let $\gamma : \Gamma \longmapsto \Gamma^\prime$ be a morphism between point-line incidence geometries. Suppose both geometries are not empty -- so they have either at least one point or at least one line. (In fact, since $\gamma$ is a map, it suffices to ask that $\Gamma$ is not empty.) Then $\gamma$ is a {\em{cover}} if $\gamma$ is locally a bijection, that is, if $x$ is any point of $\Gamma$, $\gamma$ induces a bijection between the lines incident with $x$ and the lines incident with $\gamma(x)$, and if $L$ is a line of $\Gamma$, it induces a bijection between the points incident with $L$ and the points incident with $\gamma(L)$. Sometimes we also say that $(\Gamma, \gamma)$ is a {\em{cover}} (of $\Gamma^\prime$), or even that $\Gamma$ is a {\em{cover}} of $\Gamma^\prime$ if the covering map is clear.

If $\Gamma^\prime$ is connected, then a cover $\gamma : \Gamma \longmapsto \Gamma^\prime$ is necessarily surjective; see \cite{part1}.

If $\gamma : \Gamma \longmapsto \Gamma^\prime$ is a cover, and each fiber (of lines {\em{and}} points) has constant size $\theta$, we say that $\gamma$, or $(\Gamma, \gamma)$, or $\Gamma$ is a {\em {$\theta$-fold cover}} (of $\Gamma^\prime$) or simply a {\em {$\theta$-cover}} (of $\Gamma^\prime$). We also say that $\Gamma^\prime$ is {\em {$\theta$-covered}} by $\Gamma$.

\subsection{Factorization of morphisms I -- Lower decomposition}\label{intro3}

The following theorem was proved in \cite{part1}.

\begin{theorem}[Lower decomposition]
\label{t1}
Let $\mathcal{S}$ be a thick generalized quadrangle, and let $\mathcal{S^\prime}$ be a thick full subquadrangle which is a geometric hyperplane of $\mathcal{S}$. Define $\mathcal{A}$ and $\mathcal{E}$ as above, and let $\pi : \mathcal{A} \longmapsto \mathcal{E}$ be the natural projection. Then any cover $\gamma : \mathcal{A} \longmapsto \mathcal{E}$ factorizes as $\gamma = \alpha \circ \pi$, with $\alpha$ an automorphism of $\mathcal{E}$.
\end{theorem}

\subsection{The automorphism $\overline \alpha$}\label{intro4}
Use the notation of the previous section. The following is taken from \cite{part1}.

For each line $L$ of $\mathcal{E}$ the lines of $\gamma^{-1}(L)$ all contain the same point $u$ of $\mathcal{S^\prime}$. If the ovoids of the rosette $L$ all share the point $u^\prime$, then for any other rosette $M \not= L$ whose ovoids share $u^\prime$, the lines of $\gamma^{-1}(M)$ contain $u$. If $u^\prime = \zeta(u)$, and we let $u$ vary, we obtain a permutation
\begin{equation}
\zeta : w \longmapsto w^\prime
\end{equation}
of the points of $\mathcal{S}^\prime$. This map $\zeta$ defines an automorphism $\overline \alpha$ of $\mathcal{S^\prime}$.

\subsection{Factorization of morphisms II -- Higher decomposition}\label{intro5}

We keep using the notation of the previous sections.

We say that $\pi : \mathcal{A} \longmapsto \mathcal{E}$ has the {\em {higher decomposition property}} if any cover $\gamma : \mathcal{A} \longmapsto \mathcal{E}$ factorizes as $\gamma = \pi \circ \widetilde{\alpha}$ for some $\widetilde{\alpha}  \in \Aut(\mathcal{S})_\mathcal{S^\prime}$. We also say that $\gamma : \mathcal{A} \longmapsto \mathcal{E}$ has the {\em {higher decomposition property}} 
if it factorizes as above. Such an automorphism $\widetilde\alpha$ extends the automorphism $\overline \alpha$ of $\mathcal{S^\prime}$.

\begin{theorem}\label{t2}
The cover $\pi : \mathcal{A} \longmapsto \mathcal{E}$ has the higher decomposition property if and only if any automorphism of $\mathcal{E}$ is induced by an automorphism of $\mathcal{S}$.
\end{theorem}

Note that any automorphism of $\mathcal{E}$ is induced by an automorphism of $\mathcal{S^\prime}$. If every automorphism of $\mathcal{S^\prime}$ extends to an automorphism of $\mathcal{S}$, then every automorphism of $\mathcal{S^\prime}$ induces an automorphism of $\mathcal{E}$ in a faithful manner. In such a case higher decomposition is possible.

{\bf Examples}
\begin{itemize}
\item[{\rm (1)}] Consider a subquadrangle $\mathcal{S^\prime} \cong \mQ(4, q)$ of the generalized quadrangle $\mQ(5, q) = \mathcal{S}$. Here the geometry $\mathcal{E}$ is 2-covered by $\mathcal{S} \setminus \mathcal{S^\prime} = \mathcal{A}$. As every automorphism of $\mathcal{S^\prime}$ extends to an automorphism of $\mathcal{S}$, higher decomposition is possible.
\item[{\rm (2)}] Consider a subquadrangle $\mathcal{S^\prime} \cong \mH(3, q^2)$ of the generalized quadrangle $\mH(4, q^2) = \mathcal{S}$. Here the geometry $\mathcal{E}$ is $(q + 1)$-{\em{covered}} by $\mathcal{S} \setminus \mathcal{S^\prime} = \mathcal{A}$. As every automorphism of $\mathcal{S^\prime}$ extends to an automorphism of $\mathcal{S}$, higher decomposition is possible.
\item[{\rm (3)}] Consider a subquadrangle $\mathcal{S^\prime} \cong \mQ(3, q)$ of the generalized quadrangle $\mQ(4, q) = \mathcal{S}$. Here the geometry $\mathcal{E}$ is $\theta$-covered by $\mathcal{S} \setminus \mathcal{S^\prime} = \mathcal{A}$, with $\theta = 2$ for $q$ odd and $\theta = 1$ for $q$ even. For $\theta = 1$ it is clear that higher decomposition is possible. In the odd case not every automorphism of $\mathcal{S^\prime}$ extends to an automorphism of $\mathcal{S}$. But by Theorem 10.3 of \cite {part1} also here higher decomposition is possible.
\item[{\rm (4)}] It is well known that the generalized quadrangle of order $(q, q^2)$ arising from a Kantor-Knuth flock \cite{TGQ} has subquadrangles of order $q$ isomorphic to $\mQ(4, q)$. For some of these subquadrangles the geometry $\mathcal{E}$ is 2-covered by $\mathcal{A}$; for other subquadrangles $\mathcal{E}$ is 1-covered by $\mathcal{A}$. In Section \ref{s2} higher decomposition in the 2-covered case will be considered.
\end{itemize}

\subsection{Semi partial geometries and covers}\label{intro6}

As will be briefly explained in the next theorem, some $\theta$-covers give rise to semi partial geometries. 

\begin{theorem}[\cite{part1}]
\label{t3}
Let $\mathcal{S}$ be a generalized quadrangle of order $(s, t)$ and let $\mathcal{S^\prime}$ be a subquadrangle of order $(s, t^\prime)$. If every subtended ovoid $\mO_x$ of $\mathcal{S^\prime}$ is $\theta$-subtended, $\theta > 1$, with $t^\prime \not= 1, t = st^\prime$ and $(\theta - 1)t = s^2$, then the geometry $\mathcal{E}$ is a semi partial geometry with parameters
\begin{equation}
s^\ast = s -1, t^\ast = t, \alpha^\ast = \theta, \mu^\ast = \theta(t - t^\prime).
\end{equation}
\end{theorem}

\subsection*{Examples}

\begin{itemize}
\item[{\rm (1)}] Considering Examples (1) and (4), with $\theta = 2$, of section \ref{intro5}, there arises semi partial geometries with parameters
\begin{equation}
s^\ast = q - 1, t^\ast = q^2, \alpha^\ast = 2, \mu^\ast = 2q(q -1).
\end{equation}
\item[{\rm (2)}] Considering Example (2), with $\theta = q + 1$, of section \ref{intro5}, there arises a semi partial geometry with parameters
\begin{equation}
s^\ast = q^2 -1, t^\ast = q^3, \alpha^\ast = q + 1, \mu^\ast = q(q + 1)(q^2 - 1).
\end{equation}
\end{itemize}

\begin{remark}\label{r1}
{\rm For the geometries described in the previous examples, we also refer to \cite{Brown} and \cite{HT2}}.
\end{remark}

\medskip
\section{Basic properties of Kantor-Knuth quadrangles and Kantor-Knuth ovoids}
\label{recap}

Let $\F_q$ be a finite field with $q$ odd, and let $\sigma \in \Aut(\F_q)$. Let $m$ be a given nonsquare in $\F_q$. Then one constructs a generalized quadrangle $\Gamma(q,\sigma)$ of order $(q,q^2)$, called {\em Kantor-Knuth generalized quadrangle}, as in \cite[\S 4.5 and \S 5.7]{TGQ}. As our notation suggests, this  construction is independent of the choice of $m$. On the other hand, two Kantor-Knuth GQs $\Gamma(q,\sigma)$ and $\Gamma(q,\sigma')$ are isomorphic if and only if $\sigma = \sigma'$. (See \cite{TGQ} and the references therein.) \\

In this section, $\Gamma = \Gamma(q,\sigma)$ is a nonclassical Kantor-Knuth GQ of order $(q,q^2)$, where $q = p^h$ and $p$ is an odd prime. 

\subsection{The orbits $\Omega_1$ and $\Omega_2$}
\label{recapor}

Let $\Omega$ be the set of all $q^3 + q^2$ subGQs of order $q$ of $\Gamma$; all these subGQs are isomorphic to $\mQ(4,q)$ (see chapter 5 of \cite{TGQ}, especially Theorem 5.1.9). 
In its action of $\Aut(\Gamma)$ on $\Omega$, there are two orbits $\Omega_1$ and $\Omega_2$, the first one having size $2q^2$ and 
the latter having size $(q - 1)q^2$. 

Two more  facts will be used without further notice:
\begin{itemize}
\item
$\Aut(\Gamma)$ fixes some line which we denote by $[\infty]$, and which is contained in each element of $\Omega$ {(by \cite{Glasg}, there is a line $[\infty]$ each 
point of which is a translation point, and if this line would not be fixed by $\Aut(\Gamma)$, then each point of $\Gamma$ would be a translation point; in that case, it is well known that $\Gamma$ is classical \cite{KTHVM})};
\item
$\Aut(\Gamma)$ acts transitively on the $q^4$ subGQs of order $(q,1)$ containing $[\infty]$. This can easily be derived from the fact 
that each point incident with $[\infty]$ is a translation point (see \cite{Glasg}), or equivalently that each line in $[\infty]^{\perp}$ is an axis 
of symmetry; {$\Aut(\Gamma)$ acts transitively on the line pairs $(U,V)$ for which $U, V \in [\infty]^{\perp}$ and $U \not\sim V$}.
\end{itemize}

The latter fact will be especially useful in several calculations involving orbits which we will carry out below. Other properties will be mentioned in due course.

\subsection{The automorphism group}

The size of $\Aut(\Gamma)$ is $(q + 1)(q - 1)^2q^6 \cdot \delta h$, where $\delta = 4$ if $\sigma^2 = 1$ and $\sigma \ne 1$, and $\delta = 2$ if $\sigma^2 \ne 1$ \cite{Paynecoll}.

\subsection{Subtended ovoids}

Let $\mQ \in \Omega_1$ and let $x$ be a point of $\Gamma$ which is not a point of $\mQ$; then each ovoid $\mO_x$ in $\mQ$ is doubly subtended \cite[Theorem 11.1]{part1}; 
as such, there arise $(q + 1)q^2(q - 1)/2$ subtended ovoids. This fact is not true for subGQs coming from the orbit $\Omega_2$; 
there, each subtended ovoid is $1$-subtended | see \cite[Theorem 11.1]{part1} | and so we obtain $(q + 1)q^2(q - 1)$ such ovoids. 

Call an ovoid arising from a point $x$ and a GQ in $\Omega$  a ``Kantor-Knuth ovoid.'' From the facts below, this notion is well defined, independent of the 
choice of  GQ-orbit:

\begin{itemize}
\item[FA]
let $(\mQ,x)$ and $(\mQ',x')$ be subGQ-point pairs of $\Gamma$ such that the subGQ is an element of $\Omega_1$, and the point is not contained in the subGQ; then by 
 Thas \cite{Glasg}, there is an automorphism $\alpha$ of $\Gamma$ which sends $\mQ$ to $\mQ'$ and $x$ to $x'$; in particular, 
$\alpha$ sends $\mO_x$ (in $\mQ$) to $\mO'_{x'}$ (in $\mQ'$); the same statement is true for the orbit $\Omega_2$ (we will also obtain this result in \S \ref{intr});
\item[FB]
by the previous point, by putting $\mQ = \mQ'$ (with $\mQ \in \Omega$), we conclude that $\Aut(\Gamma)_{\mQ}$ acts transitively on the points of $\Gamma$ outside $\mQ$.
\end{itemize}

By these facts, each subtended ovoid obtained from a subGQ in $\Omega$, and an exterior point (to the subGQ), merits the property of being named ``Kantor-Knuth ovoid.''

\medskip
\section{Synopsis of translation ovoids}
\label{tranov}

An ovoid $\mO$ of $\mQ(4,q)$ is a {\em translation ovoid} with respect to the point $\omega \in \mO$  if there is a point $\omega \in \mO$ such that there is a subgroup $T$ of 
$\Aut(\mQ(4,q))_\mO$ which fixes $\omega$ linewise, and which acts sharply transitively on the points of $\mO \setminus \{ \omega \}$. Equivalent 
definitions can be found in \cite{bloem}. 

If $\Aut(\mQ(4,q))_\mO$ does not fix $\omega$, it acts $2$-transitively on the points of $\mO$, and then it is an elliptic quadric in $\mQ(4,q)$ in case $q$ is odd (see e. g. \cite{2tranov}).  

If $\Gamma$ is a semifield flock TGQ of order $(q,q^2)$, and $\mQ$ is a $\mQ(4,q)$-subGQ of $\Gamma$, then it can be shown that 
any subtended ovoid in $\mQ$ is a translation ovoid. Easy way to see this: if $\Gamma$ is classical, then such ovoids 
are elliptic quadrics; if $\Gamma$ is not classical, then there is a unique line $[\infty]$ of translation points \cite{Glasg}, and this is a line of $\mQ$. Let $e$ be any 
point of $\Gamma$ not contained in $\mQ$, and let $e'$ be the point incident with $[\infty]$ which is collinear with $e$. Let $T(e')$ be the translation group of $\Gamma$ 
with translation point $e'$. For any two points $x, y$ in $\mO_e$, supposed to be not incident with $[\infty]$, there is precisely one element in $T(e')$
which maps $x$ to $y$. It is obvious that this element stabilizes $\mQ$, and so also $\mO_e$. It is now easy to see that $\mO_e$ is a translation ovoid 
relative to the point $e'$, and that $T(e')_e$ acts (sharply) transitively on $\mO \setminus \{ e' \}$.

 \section{Decomposition laws for Kantor-Knuth generalized quadrangles}
 \label{s2}

Consider a thick generalized quadrangle $\mathcal{S}$ of order $(s, s^2)$ satisfying the following conditions:
\begin{itemize}
\item[{\rm (i)}] $\mathcal{S}$ has a line $L$ of translation points;
\item[{\rm (ii)}] $\mathcal{S}$ has a subquadrangle $\mathcal{S^\prime}$ of order $s$ containing $L$;
\item[{\rm (iii)}] every subtended ovoid in $\mathcal{S^\prime}$ is doubly subtended.
\end{itemize}

By  section \ref{intro6} we know that the corresponding geometry $\mathcal{E}$ is a semi partial geometry with parameters
\begin{equation}
s^\ast = s - 1, t^\ast = s^2, \alpha^\ast = 2, \mu^\ast = 2s(s - 1).
\end{equation}

\begin{theorem}\label{t4}
Let $\gamma: \mathcal{A} \longmapsto \mathcal{E}$ be a 2-cover of $\mathcal{E}$. Then, by section \ref{intro4}, $\gamma$ induces an automorphism $\overline \alpha$ of $\mathcal{S^\prime}$. If the automorphism $\overline \alpha$ fixes the line $L$, then there exist exactly two automorphisms $\widetilde \alpha$ of $\mathcal{S}$ which extend the automorphism $\overline \alpha$ of $\mathcal{S^\prime}$.
\end{theorem}

{\em Proof}.\quad
As by Chapter 8 of \cite{PT2} all the lines of $\mathcal{S^\prime}$ concurrent with $L$ are regular, it follows that all the lines of $\mathcal{S^\prime}$ are regular; see 1.3.6 of \cite{PT2}. Hence the generalized quadrangle $\mathcal{S^\prime}$ is isomorphic to the classical generalized quadrangle $\mQ(4, s)$; see 5.2.1 of \cite{PT2}.

We will construct an automorphism $\widetilde\alpha$ of $\mathcal{S}$ such that $\overline\alpha$ is the restriction of $\widetilde\alpha$ to $\mathcal{S^\prime}$.

Choose a point $z$ of $\mathcal{A}$ and choose a line $N$ of $\mathcal{S}$ incident with $z$; the point of $\mathcal{S^\prime}$ incident with $N$ is denoted by $x^\ast$. Assume that $x^\ast$ is not on $L$. The point $z$ covers an ovoid $\mO(z)$ of $\mathcal{S^\prime}$ and $\mO(z)$ is subtended by the points $z^\prime, z^{\prime\prime}$ of $\mathcal{S}$. Let $z^{\widetilde\alpha} = z^\prime$. Then the line $N^{\widetilde\alpha}$ is the line $N^\prime = z^{\prime}x$ with $x = x^{\ast \widetilde \alpha} = x^{\ast \overline \alpha}$. For any line $T$ of $\mathcal{S^\prime}$, we define $T^{\widetilde\alpha} = T^{\overline\alpha}$.

Now let $y$ be a point of $zx^\ast$, with $x^\ast \not= y \not= z$. The point $y$ covers an ovoid $\mO(y)$ of $\mathcal{S^\prime}$ and $\mO(y)$ is subtended by a unique point $y^\prime$ of $N^\prime$; let $y^\prime = y^{\widetilde \alpha}$. Next, consider a point $v$ of $\mathcal{A}$, which is not collinear with $x^\ast$. Let $V$ be the line incident with $v$ and concurrent with $N$, let $w^\ast$ be the common point of $V$ and $\mathcal{S^\prime}$, and let $a$ be the common point of $N$ and $V$. The point $v$ covers an ovoid $\mO(v)$ of $\mathcal{S^\prime}$ and $\mO(v)$ is subtended by a point $v^\prime$ of the line $V^\prime = wa^{\widetilde \alpha}$ (note that this is indeed a line), with $w^{\ast \widetilde \alpha} = w.$  Let $V^\prime = V^{\widetilde \alpha}$ and $v^\prime = v^{\widetilde \alpha}$. Then $\widetilde \alpha$ is defined for all points of $\mathcal{S}$ not collinear with $x^\ast$, for all points of the line $N$, for all points of $\mathcal{S^\prime}$, for the line $N$, for all lines concurrent with $N$ but not containing $x^\ast$, and for all lines of $\mathcal{S^\prime}$; for all these points and lines $\widetilde \alpha$ preserves incidence.

Each line $U$ concurrent with $L$, with $L$ consisting of translation points of $\mathcal{S}$, is regular. Let $U$ be a line of $\mathcal{S}$ concurrent with $L$, not concurrent with $N$, and not contained in $\mathcal{S^\prime}$. Let $u^\ast$ be the common point of $U$ and $L$. The line $U$ of $\mathcal{A}$ covers a rosette consisting of $s$ ovoids $\mO(u_1), \mO(u_2), \ldots, \mO(u_s)$ of $\mathcal{S^\prime}$, with $u_1, u_2, \ldots, u_s$ the points of $\mathcal{A}$ incident with $U$. These ovoids have as common point $u = u^{\ast \widetilde \alpha}$, with $u$ on $L$. The ovoid $\mO(u_i)$ is subtended by points $u^1_i, u^2_i$ of lines $U_1, U_2$ incident with $u$, with $i = 1, 2, \ldots, s.$ The points $u_i$ which are not collinear with $x^\ast$ are mapped by $\widetilde \alpha$ onto points $u_i^{\widetilde \alpha}$, where $u_i^{\widetilde \alpha}$ is incident with a line $U_1, U_2$. Assume, by way of contradiction, that for at least one such $u_k$ we have that $u_k^{\widetilde \alpha}$  is incident with $U_1$, and that for at least one such $u_l$ we have that $u_l^{\widetilde \alpha}$ is incident with $U_2$. Let $\overline{U}_k$ be the line incident with $u_k$ and concurrent with $N$, and let $\overline{U}_l$ be the line incident with $u_l$ and concurrent with $N$. Let $\overline{u}_k$ be the common point of $\overline{U}_k$ and $\mathcal{S^\prime}$, and let $\overline{u}_l$ be the common point of $\overline{U}_l$ and $\mathcal{S^\prime}$. Consider also the elements $\overline{U}_k^{\widetilde \alpha}, \overline{U}_l^{\widetilde \alpha}, \overline{u}_k^{\widetilde \alpha}, \overline{u}_l^{\widetilde \alpha}$; the lines $\overline{U}_k^{\widetilde \alpha}, \overline{U}_l^{\widetilde \alpha}$ are concurrent with $N^{\widetilde \alpha}$.

First assume that $x^\ast$ and $u^\ast$ are collinear. By the regularity of $U$ the line $A$ incident with $\overline{u}_l$ and concurrent with ${x^\ast}{u^\ast}$, intersects the line
$\overline{U}_k$. As $A$ is a line of $\mathcal{S^\prime}$, it contains the point $\overline{u}_k$. Hence $\overline{u}_l \sim \overline{u}_k$. It follows that $\overline{u}_l^{\widetilde \alpha} \sim \overline{u}_k^{\widetilde \alpha}$, and that $xu$ is concurrent with $\overline{u}_l^{\widetilde \alpha}\overline{u}_k^{\widetilde \alpha}$. The lines $xu, U_1, U_2$ are regular, hence the lines $xu, \overline{u}_l^{\widetilde \alpha}\overline{u}_k^{\widetilde \alpha}, N^{\widetilde \alpha}, \overline{U}_k^{\widetilde \alpha}, \overline{U}_l^{\widetilde \alpha}, u_k^{\widetilde \alpha}u = U_1, u_l^{\widetilde \alpha}u = U_2$ belong to a common grid. Consequently $u_k^{\widetilde \alpha}u = U_1 \sim \overline{U}_l^{\widetilde \alpha}$ and $u_l^{\widetilde \alpha}u = U_2 \sim \overline{U}_k^{\widetilde \alpha}$. So there arises a triangle with sides $U_1, U_2, \overline{U}_l^{\widetilde \alpha}$, a contradiction. It follows that we may assume that $x^\ast \not \sim u^\ast$.

There exists a translation $\theta$ of $\mathcal{S}$ with translation point on $L$, which maps $u^\ast$ onto $u$, and which maps $u_k$ onto $u_k^{\widetilde \alpha}$. Then ${\widetilde \alpha}^{-1}{\theta}$ fixes the point $u$, the point $u_k^{\widetilde \alpha}$ and induces an isomorphism from $\mathcal{S^\prime}$ onto $\mathcal{S^\prime}^\theta$. The point $u$ is a translation point of the generalized quadrangle $\mathcal{S}$. We may put $\mathcal{S} = \hT(\mO)$, with $\mO$ a generalized ovoid of $\PG(4n -1, q), s = q^n$, and with $u$ the point $\infty$ of $\hT(\mO)$. The subquadrangle $\mathcal{S^\prime}$ is a $\hT(\mO^\prime)$, with $\mO^\prime \subset \mO$ a generalized conic containing $L$ in $\PG(3n - 1, q)$ (see Chapter 7 of \cite{TGQ}); the subquadrangle $\mathcal{S^\prime}^\theta$ is a $\hT(\mO^{\prime\prime}) = T(\mO^\prime)^\theta$, with $\mO^{\prime\prime} \subset \mO$ a generalized conic containing $L$. The lines $U_1 = U^\theta, U_2$ are elements of $\mO \setminus \mO^\prime$. The line $N^\theta$ is an $n$-dimensional space containing an element $\widetilde{N}^\theta \in \mO \setminus \mO^\prime, \widetilde{N}^\theta \not= U_1, U_2$, the line $N^{\widetilde \alpha}$ is an $n$-dimensional space containing an element $\widetilde{N}^{\widetilde \alpha} \in \mO \setminus \mO^\prime, \widetilde{N}^{\widetilde \alpha} \not= U_1, U_2$. We have $N^\theta \cap \mathcal{S^\prime}^\theta = \lbrace x^{\ast \theta} \rbrace$ and $N^{\widetilde \alpha} \cap \mathcal{S^\prime} = \lbrace x \rbrace$ (the point $x$ belongs to the space $\PG(3n, q)$ containing $\mathcal{S^\prime}$, and the point $x^{\ast\theta}$ belongs to the space $\PG(3n, q)^\theta$ containing $\mathcal{S^{\prime\theta}}$). The lines $\overline{U}_k^\theta, \overline{U}_l^\theta$ are $n$-spaces containing $U^\theta = U_1 \in \mO$, and they intersect $\mathcal{S^\prime}^\theta$ in points $\overline{u}_k^\theta, \overline{u}_l^\theta$. If $\langle N^\theta, U^\theta \rangle \cap \PG(3n - 1, q) = \pi_1$, with $\pi_1$ $(n - 1)$-dimensional, then $\overline{u}_k^\theta, \overline{u}_l^\theta, x^{\ast \theta}$ are contained in an $n$-space $\widetilde \pi_1$ containing $\pi_1$ (as $\overline{U}_k^\theta$ and $\overline{U}_l^\theta$ intersect $N^\theta$). The lines $\overline{U}_k^{\widetilde \alpha}, \overline{U}_l^{\widetilde \alpha}$ are $n$-spaces containing $U_1, U_2 \in \mO$, and they intersect $\PG(3n, q)$ in points $\overline{u}_k^{\widetilde \alpha}, \overline{u}_l^{\widetilde \alpha}$. If $\langle N^{\widetilde \alpha}, U_1 \rangle \cap \PG(3n - 1, q) = \zeta_1$ and $\langle N^{\widetilde \alpha}, U_2 \rangle \cap \PG(3n - 1, q) = \zeta_2$, then $\zeta_1 \cap \zeta_2 = \emptyset$, $\langle \zeta_1, \overline{u}_k^{\widetilde \alpha} \rangle =\widetilde \zeta_1, \langle \zeta_2, \overline{u}_l^{\widetilde \alpha} \rangle = \widetilde \zeta_2$ and $\widetilde \zeta_1 \cap \widetilde \zeta_2 = \lbrace x \rbrace$ (as $\overline{U}_k^{\widetilde\alpha}$ and $\overline{U}_l^{\widetilde\alpha}$ intersect $N^{\widetilde\alpha}$). Considering all the points $u_r$ of $\mathcal{S} \setminus \mathcal{S^\prime}$ on $U$, which are not collinear with $x^\ast$, there arise in $\hT(\mO)$ $s - 1$ points $\overline{u}_r^\theta$ of $\widetilde \pi_1 \setminus (\pi_1 \cup \lbrace x^{\ast \theta} \rbrace)$. Then the isomorphism $\theta^{-1} \widetilde \alpha$ of $\mathcal{S^\prime}$ onto $\mathcal{S^\prime}^\theta$ maps these points onto $s - 1$ points $\overline{u}_r^{\widetilde\alpha}$ of $(\widetilde \zeta_1 \cup \widetilde \zeta_2) \setminus (\zeta_1 \cup \zeta_2 \cup \lbrace x \rbrace)$; also, these points are not all in either $\widetilde \zeta_1$ or $\widetilde \zeta_2$. Notice also that $(x^{\ast \theta})^{\theta^{-1}\widetilde \alpha} = x$. Clearly, this yields a contradiction. Hence the $s$  points of $U$ not collinear with $x^\ast$ are mapped by $\widetilde \alpha$ onto $s$ points of one line $U_i, i \in \lbrace 1, 2 \rbrace$; this line is denoted by $U^{\widetilde \alpha}$. It follows that $\widetilde \alpha$ is defined for all lines $U$ concurrent with $L$.

Let $U$ be a line concurrent with $L$, but not concurrent with $N$ (possibly $U = L$). Consider the grid $\mathcal{G}$ containing $U$ and $N$, and let $V$ be a line of $\mathcal{G}$ not in $\mathcal{S^\prime}$ and not concurrent with $N$ or $U$. Let $x_1, x_2, \ldots, x_s$ be the points of $V$ not collinear with $x^\ast$, and let $X_1, X_2, \ldots, X_s$ be the lines of $\mathcal{G}$ concurrent with $N$, with $x_i$ incident with $X_i, i = 1, 2, \ldots, s$. The lines $N^{\widetilde \alpha}, X_1^{\widetilde \alpha}, X_2^{\widetilde \alpha}, \ldots, X_s^{\widetilde \alpha}, U^{\widetilde \alpha}$ belong to a grid $\mathcal{G}^{\widetilde \alpha}$, and $x_i^{\widetilde \alpha}$ is incident with $X_i^{\widetilde \alpha}, i = 1, 2, \ldots, s.$ The points $x_1^{\widetilde \alpha}, x_2^{\widetilde \alpha}, \ldots, x_s^{\widetilde \alpha}$ of the grid $\mathcal{G}^{\widetilde \alpha}$ belong to two concurrent lines $V_1, V_2$ (whose common point belongs to $\mathcal{S^\prime}$). It easily follows that $x_1^{\widetilde \alpha}, x_2^{\widetilde \alpha}, \ldots, x_s^{\widetilde \alpha}$ are incident with one of the lines $V_1, V_2$, say $V_1$. Let us put $V_1 = V^{\widetilde \alpha}$.

Let $\mathcal{T}$ be the grid containing $N$ and $L$, and let $Y$ be a line containing a point $y$ of $\mathcal{T}$, with $x^\ast$ not incident with $Y$. We will show that the line $Y^{\widetilde \alpha}$ is uniquely defined. We may assume that $Y \not \sim N$, that $Y \not \sim L$, and that $Y \not \in \mathcal{T}$ (putting $U = L$ in the previous paragraph, the images of all lines of $\mathcal{T}$ are well defined). Consider the grid $\mathcal{T^\prime}$ defined by $N$ and $Y$; notice that the pair $\lbrace N, Y \rbrace$ is regular as there is a (regular) line $U^\prime$ (containing $y$) intersecting $N, Y, L$. Let $u^\prime$ be the common point of $U^{\prime}$ and $L$. The grid $\mathcal{T^\prime}$ has a line $U^{\prime\prime}$ containing $u^{\prime}$, with $U^{\prime\prime} \not \sim Y$. By the foregoing paragraph the line $Y^{\widetilde \alpha}$ is well defined.


Let $Y$ be a line not containing a point of $\mathcal{T}$ and not belonging to $\mathcal{S^\prime}$. Let $y$  be the common point  of $Y$ and $\mathcal{S^\prime}$, let $Z$ be the line containing $y$ and concurrent with $L$, and let $Z \cap L = \lbrace z \rbrace$. Further, let $U$ be the line containing $z$ and concurrent with $N$. The grid defined by the regular line $U$ and the line $Y$ is denoted by $\mathcal{R}$. We determine the intersection of $\mathcal{R}$ and $\mathcal{S^\prime}$. 

First assume that $x^\ast \not \sim z$. There are $s$ points $x_1, x_2, \ldots, x_s$ of $\mathcal{R}$ in $\mathcal{S^\prime}$, which are not incident with $Z$. The line $Z^\prime$ containing $x_1$ and concurrent with $Z$ belongs to $\mathcal{R}$ and $\mathcal{S^\prime}$, so contains the points $x_1, x_2, \ldots, x_s$. The line $Z^\prime$ does not contain $y$, as otherwise $\mathcal{R}$ contains at least three lines through $y$. The lines $Z, Z_1, Z_2, \ldots, Z_s$ of $\mathcal{R}$ concurrent with $U$ have well-defined images $Z^{\widetilde \alpha}, Z_1^{\widetilde \alpha}, Z_2^{\widetilde \alpha}, \ldots, Z_s^{\widetilde \alpha}$ as they have a point in common with $\mathcal{T}$. Also the regular line $U$ and the line $Z^\prime$, both belonging to $\mathcal{R}$, have well-defined images $U^{\widetilde \alpha}$ and $(Z^\prime)^{ \widetilde \alpha}$. Clearly $U^{\widetilde \alpha}$ is also regular. Hence $Z^{\widetilde \alpha}, Z_1^{\widetilde \alpha}, Z_2^{\widetilde \alpha}, \ldots, Z_s^{\widetilde \alpha}, U^{\widetilde \alpha}$ and $(Z^\prime)^{ \widetilde \alpha}$ belong to a grid $\mathcal{R^{\widetilde \alpha}}$. 
The images of the points of $Y$ are incident with lines $Y_1, Y_2$ and belong to $\mathcal{R^{\widetilde \alpha}}$; the common point of $Y_1, Y_2$ is the point $y^{\widetilde \alpha}$. It easily follows that these images are incident with just one of the lines $Y_1, Y_2$, say $Y_1$. So $Y^{\widetilde \alpha} = Y_1$ is well defined.

Next, assume that $x^\ast \sim z$. Then the intersection of $\mathcal{R}$ and $\mathcal{S^\prime}$ consists of the lines $U$ and $Z$. Let $W$ be a line of $\mathcal{S^\prime}$ through $y$, not containing $z$, and let $\mathcal{R^\prime}$ be a grid containing $W, Y$, and a point of $L$, but not containing $x^\ast$. Then all lines of $\mathcal{R^\prime}$ distinct from $Y$ have a well-defined image under $\widetilde\alpha$ and belong to a grid $\mathcal{R^{\prime \widetilde \alpha}}$. It easlly follows that $Y^{\widetilde \alpha}$ is well defined.

We have shown that $W^{\widetilde \alpha}$ is well defined for each line of $\mathcal{S}$ which does not contain $x^\ast$.

We have still to define $t^{\widetilde \alpha}$, with $t \sim x^\ast$, $t$ not incident with $N$, $t$ not in $\mathcal{S^\prime}$, and also $W^{\widetilde \alpha}$, with $x^\ast$ incident with $W$, $W \not= N$, $W$ not in $\mathcal{S^\prime}$. This is easy. Let $t$ be such a point. Let $E$ be a line containing $t$, with $E \not= tx^\ast$, and let $t^\prime$ be the point incident with $E^{\widetilde \alpha}$ and collinear with $x^{\ast \widetilde \alpha} = x$. Now let $E^\prime$ be a second line containing $t$, with $E^\prime  \not= tx^\ast$, and let $t^{\prime \prime}$ be the point incident with $E^{\prime \widetilde \alpha}$ and collinear with $x$. Assume, by way of contradiction, that $t^\prime \not= t^{\prime \prime}$. If $E^{\widetilde \alpha}$ and $E^{\prime \widetilde \alpha}$ would be incident with a common point $\widehat{t}$, then, as $t^\prime \not= t^{\prime \prime}$, the point $\widehat{t} ^{{\widetilde \alpha}^{-1}}$ would be incident with $E$ and $E^\prime$, with $\widehat{t} ^{{\widetilde \alpha}^{-1}} \not= t$, a contradiction. Let $D \sim E^{\widetilde \alpha}, D \sim E^{\prime \widetilde \alpha}$, $t^\prime$ not incident with $D$ and $t^{\prime \prime}$ not incident with $D$. Then $D^{\widetilde {\alpha}^{-1}} \sim E, D^{{\widetilde \alpha}^{-1}} \sim E^\prime$, $t$ not incident with $D^{{\widetilde \alpha}^{-1}}$. Hence there arises a triangle in $\mathcal{S}$, a contradiction. Consequently $t^\prime = t^{\prime \prime}$. Let us then put $t^\prime = t^{\widetilde \alpha}$. 

Finally, let $x^\ast$ be incident with $W, W \not= N$, $W$ not in $\mathcal{S^\prime}$. Let $t$ be incident with $W$, $t \not= x^\ast$, and let $W^\prime = xt^{\widetilde \alpha}$. Let $r$ be incident with $W$, $r \not= x^\ast, r \not= t$, and let $W^{\prime \prime} = xr^{\widetilde \alpha}$. Assume, by way of contradiction, that $W^\prime \not= W^{\prime \prime}$ so that $r^{\widetilde \alpha} \not \sim t^{\widetilde \alpha}$. Let $b \sim r^{\widetilde \alpha}, b \sim t^{\widetilde \alpha}, b \not= x$. Then $b^{{\widetilde \alpha}^{-1}} \sim r, b^{{\widetilde \alpha}^{-1}} \sim t$, with $b^{{\widetilde \alpha}^{-1}}$ not incident with $W$, clearly a contradiction. Let us then put $W^\prime = W^{\widetilde \alpha}$.

Now $\widetilde \alpha$ is defined for all points and lines of $\mathcal{S}$ and preserves incidence. Hence $\widetilde \alpha$ is an automorphism of $\mathcal {S}$ such that $\overline \alpha$ is the restriction of $\widetilde \alpha$ to $\mathcal{S^\prime}$.

If in the second paragraph of the proof the point $z^{\widetilde \alpha}$ is chosen to be $z^{\prime \prime}$, then the second such automorphism $\widetilde \alpha$ is obtained. \eop \\

Before proceeding, recall the next theorem.

\begin{theorem}[K. Thas \cite{Stab, TGQ}]
\label{t5}
Let $\mathcal{S}$ be a thick generalized quadrangle of order $(s, s^2)$ having a line $L$ of translation points and a subquadrangle $\mathcal{S^\prime}$ of order $s$ containing $L$, with the property that every subtended ovoid in $\mathcal{S^\prime}$ is doubly subtended. Then, for $s$ even $\mathcal{S}$ is classical and for $s$ odd $\mathcal{S}$ is a Kantor-Knuth generalized quadrangle. 
\end{theorem}

\begin{remark}\label{r2}\cite{TGQ}
{\rm If $\mathcal{S}$ is a nonclassical Kantor-Knuth generalized quadrangle of order $(s, s^2)$, then any subquadrangle of order $s$ contains the unique line $L$ of translation points.}
\end{remark}

\begin{theorem}\label{t6}
Let $\mathcal{S}$ be  a nonclassical Kantor-Knuth generalized quadrangle, let $L$ be the line of translation points, let $\mathcal{S^\prime}$ be a subquadrangle containing $L$ and assume that every subtended ovoid in $\mathcal{S^\prime}$ is doubly subtended. If $\gamma: \mathcal{A} \longmapsto \mathcal{E}$ is a cover of $\mathcal{E}$ with $\mathcal{A} = \mathcal{S} \setminus \mathcal{S}^\prime$, then the  corresponding automorphism $\overline \alpha$ of $\mathcal{S^\prime}$ fixes the line $L$; also each automorphism of $\mathcal{E}$ fixes $L$.
\end{theorem}

{\em Proof}.\quad
The points of the geometry $\mathcal{E}$ are subtended ovoids of the subquadrangle $\mathcal{S^\prime} \cong \mQ(4, s)$. Each such ovoid is a translation ovoid relative to some point on $L$.
It follows that the automorphism $\overline \alpha$ of $\mathcal{S^\prime}$ induced by $\gamma$ fixes the line $L$. \eop \\

\begin{theorem}\label{t7}
Let $\mathcal{S}$ and $\mathcal{S^\prime}$ satisfy (i), (ii), (iii) as in Theorem \ref{t4}. Then the cover $\pi : \mathcal{A} \longmapsto \mathcal{E}$ has the higher decomposition property.
\end{theorem}

{\em Proof}.\quad
If $\mathcal{S}$ is classical, then by \cite{part1} $\pi$ has the higher decomposition property. Hence suppose that $\mathcal{S}$ is not classical, so that by Theorem \ref{t5} $s$ is odd and $\mathcal{S}$ is a Kantor-Knuth generalized quadrangle. Let $\gamma : \mathcal{A} \longmapsto \mathcal{E}$ be any cover of $\mathcal{E}$. Then by Theorem \ref{t6} the corresponding automorphism $\overline \alpha$ of $\mathcal{S^\prime}$ fixes the line $L$ of translation points. By Theorem \ref{t4} the automorphism $\overline \alpha$ can be extended in two ways to an automorphism $\widetilde \alpha$ of $\mathcal{S}$ for which $\gamma = \pi \circ {\widetilde\alpha}$. Consequently $\pi$ has the higher decomposition property. \eop \\

\begin{corollary}\label{c1}
Each automorphism $\alpha$ of the semi partial geometry $\mathcal{E}$ arising from $\pi$ is induced by exactly two automorphisms of the geometry $\mathcal{A}$.\eop 
\end{corollary}

\begin{theorem}\label{t8}
Let $\sigma \not= 1$ be the companion automorphism of the nonclassical Kantor-Knuth generalized qudrangle $\mathcal{S}$. Let $\mathcal{S^\prime}$ be a subquadrangle of order $s$ with $\mathcal{S^\prime} \in \Omega_1$.  
\begin{itemize}
\item
If $\sigma^2 = 1$, then $\Aut(\mathcal{E}) = \Aut(\mathcal{S^\prime})_L$ and $\vert\Aut(\mathcal{A})\vert = 2\vert\Aut(\mathcal{E})\vert$.
\item
If $\sigma^2 \not= 1$, then $2\vert\Aut(\mathcal{E})\vert = \vert\Aut(\mathcal{S^\prime})_L\vert$ and $\vert\Aut(\mathcal{A})\vert = 2\vert\Aut(\mathcal{E})\vert$.
\end{itemize}
\end{theorem}
{\em Proof}.\quad
By \cite{PW} $ \vert\Aut(\mathcal{S})\vert = (s + 1)(s - 1)^{2}s^{6}h\delta$, with $s = p^h$ and $p$ prime, $\delta = 2$ for $\sigma^2 \not= 1$ and $\delta = 4$ for $\sigma^2 = 1$. The orbit $\Omega_1$ of subquadrangles $\mathcal{S^\prime}$ of order $s$ with 2-subtended ovoids has size $2s^2$; see \cite{part1}. Hence $\vert\Aut(\mathcal{A})\vert = \vert\Aut(\mathcal{S})\vert/2s^2$. So $\vert\Aut(\mathcal{A})\vert = (s - 1)^{2}(s + 1)s^{4}h\delta / 2$. By Corollary \ref{c1} we have $2\vert\Aut(\mathcal{E})\vert = \vert\Aut (\mathcal{A})\vert$, so $\vert\Aut(\mathcal{E})\vert = (s - 1)^{2}(s + 1)s^{4}h\delta/4$.

Also $\vert\Aut(\mathcal{S^\prime})_L\vert = hs^{4}(s - 1)(s^2 - 1)$; see \S4.6.2 of \cite{POL}. Consequently, for $\sigma^2 = 1$ we have $\Aut(\mathcal{S^\prime})_L = \Aut(\mathcal{E})$, and for $\sigma^2 \not= 1$ we have $2\vert\Aut(\mathcal{E})\vert = \vert\Aut(\mathcal{S^\prime})_L\vert$. \eop \\

{\bf Alternative proof of Theorem \ref{t4} and Corollary \ref{c1} in the case of  $\sigma^2 = 1, \sigma \not=1$.} 

Assume that $\sigma^2 = 1, \sigma \not= 1$. Then $\vert\Aut(\mathcal{S})_{\mathcal{S}^\prime}\vert = \vert\Aut(\mathcal{A})\vert = 2\vert\Aut(\mathcal{S^\prime})_L\vert$, the latter equality arising by merely comparing the explicit values of $\vert \Aut(\mA) \vert = \vert \Aut(\mS)_{\mS'} \vert$ (when $\sigma^2 = 1$, $\sigma \ne 1$) and $2\vert \Aut(\mS')_L \vert$. As any subtended ovoid of $\mathcal{S^\prime}$ is doubly subtended we have that $\vert\Aut(\mathcal{A})\vert\leq 2\vert\Aut(\mathcal{E})\vert$ and, by Theorem \ref{t6}  we have $\vert\Aut(\mathcal{S^\prime})_L\vert \geq \vert\Aut(\mathcal{E})\vert$, that is $\vert\Aut(\mathcal{E})\vert \leq \vert\Aut(\mathcal{S^\prime})_L) \vert = \vert\Aut(\mathcal{A})\vert/2 \leq \vert\Aut(\mathcal{E})\vert$. Consequently $\Aut(\mathcal{S^\prime})_L = \Aut(\mathcal{E})$ and $\vert\Aut(\mathcal{A})\vert = 2\vert\Aut(\mathcal{E})\vert$. This proves Theorem \ref{t4} and Corollary \ref{c1} for $\sigma^2 = 1, \sigma \not= 1$. \eop \\

\medskip
\section{Embedded Kantor-Knuth ovoids}
\label{emb}

Let $\mQ$ and $\mO_x$ be as in section \ref{recap}.
Put $A := \Aut(\mQ) \cong \Aut(\mQ(4,q))$.  Then each element of $\mK(\mQ) := \mO_x^A$ is, by definition, also called {\em Kantor-Knuth ovoid}.

\medskip
\subsection{The embedding category}
\label{embcat}

Now consider an embedding 
\begin{equation}
\label{embeq}
\gamma:\ \mQ(4,q) \ \hookrightarrow\ \ \Gamma 
\end{equation}
which we see as an isomorphism $\gamma: \mQ(4,q) \mapsto \mQ$, where $\mQ$ is contained in $\Omega_1 \cup \Omega_2$ | say $\mQ \in \Omega_1$ w. l. o. g.  
Let $\mO(\mQ)$ be the set of subtended Kantor-Knuth ovoids in $\mQ$ from points of $\Gamma \setminus \mQ$ (and we use this notation throughout). We call $\gamma^{-1}(\mO(\mQ))$ the set of 
{\em embedded ovoids} of $\mQ(4,q)$ {\em with respect to the embedding $\gamma$}.   Now we define a category $\mathbf{\texttt{Emb}}(\mQ(4,q),\Gamma)$ in which 
objects are diagrams of the above form 
\begin{equation}
\gamma:\ \mQ(4,q) \ \rightarrow\ \ \mQ 
\end{equation}
and in which a morphism between objects $\gamma:\ \mQ(4,q) \ \rightarrow\ \ \mQ$ and  $\gamma':\ \mQ(4,q) \ \rightarrow\ \ \mQ'$ {(where $\mQ' \in \Omega$)} 
is a pair $(\alpha,\beta)$ of generalized quadrangle morphisms $\alpha:\ \mQ(4,q) \rightarrow \mQ(4,q)$ and
 $\beta:\ \mQ \rightarrow \mQ'$ such that the following diagram commutes: 

\begin{center}
\begin{tikzpicture}[>=angle 90,scale=3.5,text height=2.0ex, text depth=0.45ex]
\node (a0) at (0,3) {$\mQ$};
\node (a0') [right=of a0] {};
\node (a1) [right=of a0'] {$\mQ'$};

\node (b0) [above=of a0] {$\mQ(4,q)$};
\node (b1) [above=of a1] {$\mQ(4,q)$};

\draw[->,font=\scriptsize,thick,>=angle 45,orange]

(a0) edge node[above] {$\beta$} (a1)
(b0) edge node[above] {$\alpha$} (b1);



\draw[->,font=\scriptsize,thick]

(b1) edge node[right] {$\gamma'$} (a1);

\draw[->,font=\scriptsize,thick]
(b0) edge node[left] {$\gamma$} (a0);

\end{tikzpicture}
\end{center}


and such that $\beta(\mO(\mQ)) \subseteq \mO(\mQ')$. Without the latter property, the category would not detect the initial embeddings of the type given in  
(\ref{embeq}): take for instance one flag $(u,U)$ in the first $\mQ(4,q)$ on the top row, and 
define $\alpha$ by mapping all flags of $\mQ(4,q)$ on $(u,U)$. Define $\beta$ by mapping all flags of $\mQ$ to $(\gamma'(u),\gamma'(U))$.

Since $\gamma$ and $\gamma'$ are isomorphisms, the morphism $\beta$ also is an isomorphism if and only if $\alpha$ is. 

The identity element of the object $\gamma:\ \mQ(4,q) \ \rightarrow\ \ \mQ$ is $(\id_{\mQ(4,q)},\id_\mQ)$. \\

Later, we will show that we can relax the definition of morphism quite a bit, as follows: $(\alpha,\beta)$ is as above, but we ask that there is some 
$\mO \in \mO(\mQ)$ such that $\beta(\mO) \in \mO(\mQ')$ (instead of the entire action on $\mO(\mQ)$), and that some special line for $\mO$ is mapped to an appropriate special line for $\beta(\mO)$.

Note also 
that initially, $\alpha$ and $\beta$ not necessarily are isomorphisms. Still, we will show (in Corollary \ref{core}) that they {\em have} to be. 

First we need some more theory.

A {\em geometrical hyperplane} in a thick GQ $\mS$ of order $(s,t)$ is a subgeometry $\mS'$ with the property that each line $U$ of $\mS$ either contains 
one point of $\mS'$, and then $U$ is not a line of $\mS'$, or all its points are points of $\mS'$, and then $U$ is a line of $\mS'$. It is easy to show that there only are 
three types of geometrical hyperplanes in this setting:
\begin{itemize}
\item[\framebox{A}]
$\mS'$ is an ovoid of $\mS$ (and so $\mS'$ does not have lines);
\item[\framebox{B}]
$\mS'$ is the point-line geometry of a set $x^{\perp}$, with $x$ a point of $\mS$ (so the only lines of $\mS'$ are those lines of $\mS$ incident with $x$);
\item[\framebox{C}]
$\mS'$ is a subGQ of order $(s,t/s)$. This is the typical finite formulation of this case; in the general case, we want to say that $\mS'$ is a proper full subGQ which 
contains at least one point of every line in $\mS$.  Note that such a subGQ is maximal in the following sense: there is no proper subGQ $\widehat{\mS}$ of $\mS$ 
which properly contains $\mS'$. (Suppose by way of contradiction that such a $\widehat{\mS}$ exists. Let $u$ be a point in $\widehat{\mS} \setminus \mS'$, and 
let $U$ be any line of $\mS$ incident with $u$. Then $U$ is incident with some point $v$ of $\mS'$. As $v$ is also a point of $\widehat{\mS}$, it follows that 
$U$ is a line of $\widehat{\mS}$. So $\widehat{\mS}$ is both full and ideal, and hence coincides with $\mS$ by \cite[proposition 1.8.2]{POL}.)
\end{itemize}

In the light of the previous discussion, we now show the next result (Theorem \ref{hyps}). 
First we observe an easy property. 

\begin{observation}
\label{obs}
If $\Gamma$ is a thick GQ with an ovoid $\mO$, and $\Gamma'$ is a subGQ which contains $\mO$, then $\Gamma = \Gamma'$. 
\end{observation}

{\em Proof}.\quad 
Consider any line $U$ of $\Gamma'$. Let $u \I U$ be arbitrary in $\Gamma$, but not contained in $\mO$. Let $V \I u$, $V \ne U$. Then $V$ contains a point $\omega$ of $\mO$. Since $V$ is the unique line in $\Gamma$ incident with $\omega$ and concurrent with $U$, it follows that $u \in \Gamma'$, and so $\Gamma'$ is full. Now let $v \in \Gamma'$ and not 
contained in $\mO$. Then any line on $v$ contains a point of $\mO$, so such a line must be a line of $\Gamma'$. It follows that besides full, $\Gamma'$ is also ideal, so $\Gamma' = \Gamma$ by \cite[proposition 1.8.2]{POL}. \eop \\

We also recall an interesting result of Pasini \cite{P}, which generalizes an older result of Hughes \cite{H} (which is the case of projective planes).

\begin{theorem}[Hughes--Pasini \cite{H,P}]
\label{HP}
Let $\alpha$ be a morphism from a thick (possibly infinite) generalized $m$-gon $\mE$ to a thick (possibly infinite) generalized $m$-gon $\mE'$, with $m \geq 3$. If $\alpha$ is surjective, then 
either $\alpha$ is an isomorphism, or each element in $\mE'$ has an infinite fiber in $\mE$.
\end{theorem}

As we will see later on, the thickness condition is crucial (as the statement is not true without that assumption). \\

\begin{theorem}
\label{hyps}
Let $\Delta$ and $\Delta'$ be finite thick GQs, and let $\phi: \Delta \mapsto \Delta'$ be a GQ morphism which surjectively maps a geometrical hyperplane $\mG$ of $\Delta$ to 
a geometrical hyperplane $\mG'$ of $\Delta'$:
\begin{equation}
\phi(\mG) = \mG'.
\end{equation}
Then $\phi$ is surjective if $\mG'$ is not of type B, and also if $\mG'$ is not thin if it is of type C. 
\end{theorem}

{\em Proof}.\quad
First note that $\phi(\Delta)$ is a (possibly degenerate) GQ: if $x'$ is a point of $\phi(\Delta)$ not incident with the line $Y'$ in $\phi(\Delta)$, and 
$x \in \phi^{-1}(x')$, $Y \in \phi^{-1}(Y')$, then $x$ is not incident with $Y$, and so there is a unique line $Z$ such that $x \I Z \sim Y$. So 
there is a unique line $Z' = \phi(Z)$ in $\phi(\Delta)$ such that $x' \I' Z'\ \sim Y'$. \\

{\bf CASE 1.}\quad   
Before starting the proof, note that the image of a partial ovoid through a GQ morphism not necessarily is a partial ovoid. Example: map all flags 
of a GQ $\mS$, with partial ovoid $\mO$, to a flag of type $(v,V)$, where $V$ is a fixed line, and such that the image of $\mO$ is not a singleton. \\

We now prove the statement for ovoids $\mO = \mG'$. In that case, $\mG$ also is an ovoid. Observe that $\phi(\Delta)$ is a subGQ of 
$\Delta'$ which contains $\mG'$;  as $\mG'$ is an ovoid of the thick GQ $\Delta'$, Observation \ref{obs} shows that $\phi(\Delta) = \Delta'$, and so $\phi$ is surjective. \\

{\bf CASE 2}.\quad
Now suppose $\mG'$ is of type C. Then $\mG'$ is  thick by assumption. 
Obviously $\mG$ is also of type C. Let $x$ be any point in $\Delta \setminus \mG$. First suppose that $\phi(x) \in \mG'$ for each such point. Then $\phi$ induces an 
epimorphism $\widehat{\phi}$ from the thick GQ $\Delta$ onto the thick GQ $\mG'$. By Theorem \ref{HP}, $\widehat{\phi}$ is an isomorphism, so that 
$\Delta$ and $\mG'$ have the same parameters. As $\widehat{\phi}$ induces an isomorphism from $\mG$ to $\mG'$, also $\mG$ and $\mG'$ have the same parameters. 
Hence $\mG = \Delta$, a contradiction. \\

Finally, if for some point $x \in \Gamma \setminus \mG$, we have that $\phi(x) \not\in \mG'$, we find that by maximality of $\mG'$, $\Delta'$ and $\phi(\Delta)$ must coincide. 
\eop \\

{\bf Note}.\quad 
Now we give a counter example to Theorem \ref{hyps} for $\mG'$ a thin subGQ of the GQ $\Delta'$ of order $3$. Consider $\hT(\mO) = \Delta$ in $\PG(3,3)$ with $\mO = \{ a_1,a_2,a_3,a_4 \}$ a conic in the plane $\delta$. Let $\zeta \ne \delta$ be a plane in $\PG(3,3)$ containing $a_1$ and $a_2$. Further, let $t_i$ be the tangent line of $\mO$ at $a_i$, with $i = 1, 2, 3, 4$. By definition, the (thin) GQ $\mG = \mG'$ consists of the lines of $\Delta$ in $\zeta$, the points of $\Delta$ in $\zeta$, the point $(\infty)$ 
and the planes tangent to $\mO$ which contain $t_1$ or $t_2$.  Also, let $\Delta'$ be any GQ of order $3$ intersecting $\Delta$ in $\mG = \mG'$, or let $\Delta' = \Delta$. 
Now we construct the morphism $\phi: \Delta \mapsto \Delta'$. Let $a_1a_3 \cap a_2a_4 = \{ r \}$. Each element of $\mG$ is mapped onto itself, so $\phi$ maps surjectively $\mG$ onto $\mG'$. If $m$ is a point of $\PG(3,3)$ not in $\delta \cup \zeta$, then $\phi(m) = rm \cap \zeta$; $\phi(a_3) = a_1$ and $\phi(a_4) = a_2$; 
if $R$ is a line of $\PG(3,3)$ containing $a_i$, $i \in \{ 3, 4\}$, but not contained in $\delta$, then $\phi(R) = rR \cap \zeta$; finally, if $\eta \ne \delta$ is a plane containing $t_i$, $i \in \{ 3, 4\}$ | say $i = 3$, containing the point $\ell$ of $\PG(3,3) \setminus \delta$, then $\phi(\eta)$ is the plane $t_1\phi(\ell)$. 
  
In a forthcoming paper, we hope to handle Theorem \ref{hyps}, CASE 2 without further finiteness and/or thickness restrictions.\\ 

A {\em groupoid} is a category in which each morphism is an isomorphism. The {\em core} $\cor(\mC)$ of a category $\mC$ is the maximal groupoid contained 
in the category as a subcategory.

\begin{corollary}
\label{core}
We have that $\cor(\texttt{Emb}(\mQ(4,q),\Gamma)) = \texttt{Emb}(\mQ(4,q),\Gamma)$.
\end{corollary}

{\em Proof}.\quad
Consider a morphism $(\alpha,\beta)$. Since $\beta$ sends subtended ovoids of $\mQ$ to subtended ovoids of $\mQ'$, we have by Theorem \ref{hyps} that 
$\beta$ is surjective. Hence bijective, since we consider {\em finite} generalized quadrangles.  So $\alpha$ is an isomorphism as well, and hence  $\texttt{Emb}(\mQ(4,q),\Gamma)$ is a groupoid. \eop \\

Two natural questions arise:
\begin{itemize}
\item[\framebox{DT}]
If $\mG$ is a geometrical hyperplane of $\Delta$, and $\mG'$ is a geometrical hyperplane of $\Delta'$, can $\phi(\mG)$ have a different type than $\mG$? 
\item[\framebox{TB}]
What happens if $\mG'$ {\em is} of type B?
\end{itemize}

It turns out that the answer to the second question is, that the statement of the previous theorem in case $\mG'$ is of type $B$, is not true, and that the answer of the first question is ``yes.''

We return to question TB.
We keep using the notation of Theorem \ref{hyps} and allow infinite generalized quadrangles $\Delta$ and $\Delta'$, and 
let $\mG'$ be of type B. We also suppose that $\mG$ is of type B. Suppose $\underline{\gamma}$ is a surjective morphism from the geometry  $x^{\perp} = \mG$ to the geometry of $\underline{\gamma}(x)^{\perp} = \mG'$. Note that, if $\Delta$ has (possibly infinite) parameters $(u,v)$ and $\Delta'$ parameters $(u',v')$, then as soon as $u \geq u'$ and $v \geq v'$, such 
$\underline{\gamma}$ are easy to construct. The condition is also necessary. Now let $w$ be any point of $\Delta$ not collinear with $x$; 
define $\gamma(w) := \underline{\gamma}(x)$. Let $W$ be any line incident with $w$, and let $\mathrm{proj}_x(W) =: W'$; then define $\gamma(W) = \underline{\gamma}(W')$. 
For any point or line $X$ in $\mG$, let $\gamma(X) := \underline{\gamma}(X)$. Then $\gamma$ is a well-defined morphism $\gamma: \Delta \mapsto \Delta'$ 
which surjectively maps $\mG$ to $\mG'$. \\

We also see that Theorem \ref{HP} is not true for thin quadrangles $\mE'$ (as $\mE = \Delta$ can be taken finite). \\

Now we turn to question DT.
Now let $\mG$ be of type C, and of order $(u,v)$. Let $\mG'$ be of type B, and let the parameters of $\Delta'$ be $(u',v')$. Suppose that $u \geq u'$, $v \geq v'$, and construct 
$\gamma:\ \mG \mapsto \Delta'$ as above, where $x$ is any point of $\mG$. We will now extend $\gamma$ to a morphism $\overline{\gamma}: \Delta \mapsto \Delta'$ with image 
$\mG'$, and thus answer our first question affirmatively. 

First of all, $\overline{\gamma}$ coincides with $\gamma$ on $\mG$. 
For $w$ a point in $\Delta$ not contained in $\mG$, send $w$ to $\gamma(x)$. Now let $W$ be a line incident with $w$. If $W$ meets some line $W'$ which is incident with $x$, 
then $\overline{\gamma}(W) := \gamma(W')$. If $W$ is not incident with such a line, choose $\overline{\gamma}(W)$ arbitrarily on $\gamma(x)$. 

It is easy to see that $\overline{\gamma}$ is a morphism $\overline{\gamma}: \Delta \mapsto \Delta'$ such that 
\begin{equation}
\overline{\gamma}(\Delta) = \overline{\gamma}(\mG) = \mG'.
\end{equation}

\medskip
\subsection{The category $\widetilde{\texttt{Emb}}(\mQ(4,q),\Gamma)$}

Note that a priori, section \ref{embcat} is just one possible approach to the ``embedding category'': we could also have considered as objects the 
diagrams 
\begin{equation}
\mQ(4,q) \ \overset{\gamma}{\longrightarrow} \ \mQ \ \overset{\iota_1}{\hooklongrightarrow}\ \Gamma
\end{equation}
with $\iota_1$ the canonical embedding of $\mQ$ into $\Gamma$, 
where morphisms between objects $\mQ(4,q) \ \overset{\gamma}{\longrightarrow} \ \mQ \ \overset{\iota_1}{\hooklongrightarrow}\ \Gamma$ and $\mQ(4,q) \ \overset{\gamma'}{\longrightarrow} \ \mQ' \ \overset{\iota_2}{\hooklongrightarrow}\ \Gamma$ would be $3$-tuples $(\alpha,\beta,\delta)$ of generalized quadrangle morphisms $\alpha: \mQ(4,q) \mapsto \mQ(4,q)$,  
$\beta: \mQ \mapsto \mQ'$ 
and $\delta: \Gamma \mapsto \Gamma$ such that the diagram

\begin{center}
\begin{tikzpicture}[>=angle 90,scale=3.5,text height=2.0ex, text depth=0.45ex]
\node (a0) at (0,3) {$\mQ$};
\node (a0') [right=of a0] {};
\node (a1) [right=of a0'] {$\mQ'$};
\node (c0) [below=of a0] {$\Gamma$};
\node (c1) [below=of a1] {$\Gamma$};

\node (b0) [above=of a0] {$\mQ(4,q)$};
\node (b1) [above=of a1] {$\mQ(4,q)$};

\draw[->,font=\scriptsize,thick,>=angle 45,orange]
(b0) edge node[above] {$\alpha$} (b1)
(c0) edge node[above] {$\delta$} (c1)
(a0) edge node[above] {$\beta$} (a1);



\draw[->,font=\scriptsize,thick]

(b1) edge node[right] {$\gamma'$} (a1);

\draw[->,font=\scriptsize,thick]
(b0) edge node[left] {$\gamma$} (a0);

\draw[right hook->,dotted,font=\scriptsize,thick,red]
(a0) edge node[left] {$\iota_1$} (c0)

(a1) edge node[right] {$\iota_2$} (c1);
\end{tikzpicture}
\end{center}
commutes. Note that $\alpha$ is an isomorphism if and only if $\beta$ is, and that if $\delta$ is an isomorphism (the case which is for us the most important one), it induces an isomorphism  $\underline{\delta}: \iota_1(\mQ) \mapsto \iota_2(\mQ')$, so that $\beta$ and $\alpha$ also are isomorphisms. (If $\kappa = (\alpha,\beta,\delta)$ is an isomorphism, then there should be  a $\kappa^{-1}$ such that $\kappa^{-1} \circ \kappa = (\id_{\mQ(4,q)},\id_\mQ,\id_\Gamma)$. This can only happen if and only if 
$\alpha, \beta$ and $\delta$ are isomorphisms.)

The identity element of the object $\mQ(4,q) \ \overset{\gamma}{\longrightarrow} \ \mQ \ \overset{\iota_1}{\hooklongrightarrow}\ \Gamma$ is $(\id_{\mQ(4,q)},\id_{\mQ},\id_{\Gamma})$.

We  will call this second category $\widetilde{\texttt{Emb}}(\mQ(4,q),\Gamma)$; later on, we will compare it with $\texttt{Emb}(\mQ(4,q),\Gamma)$. \\

As we will later see, a very important natural question arises: \ul{to understand the (iso)morphisms between objects in $\texttt{Emb}(\mQ(4,q),\Gamma)$}. 







\medskip
\subsection{Special lines}

Consider a $\mQ(4,q)$-quadrangle $\mQ$, and let $\mO$ be a nonclassical Kantor-Knuth ovoid. Then $\mO$ has a special point which is fixed by $\Aut(\mO) = \Aut(\mQ)_\mO$, say $x$. 
The group $\Aut(\mO)$ acts transitively on the points of $\mO \setminus \{x\}$. Suppose that 
\begin{equation}
\gamma_i: \mQ\ \mapsto\ \mQ'_i
\end{equation}
is an isomorphism of generalized quadrangles, with $\mQ'_i$ a $\mQ(4,q)$-quadrangle in a nonclassical Kantor-Knuth quadrangle $\Gamma$ in orbit $\Omega_i$ ($i = 1,2$), such that $\gamma_i(\mO)$ is a subtended ovoid in $\mQ'_i$. Let $[\infty]$ be the special  line of $\Gamma$; then $[\infty] \I \gamma_i(x)$ and $[\infty]$ is a line of $\mQ'_i$. Define $U_i := \gamma_i^{-1}([\infty])$; it is a line incident with $x$. 
Then one can prove that the $\Aut(\mQ)_\mO$-orbit which contains $U_1$ is different from the $\Aut(\mQ)_\mO$-orbit which contains $U_2$; this will be done in the next section.
It is important to note that as such we obtain two distinct orbits in $\mQ$ of lines incident with $x$, say $U_1(\mO), U_2(\mO)$.  Here, $U_1(\mO)$ denotes the lines coming from subGQs in the orbit $\Omega_1$, and $U_2(\mO)$ corresponds to subGQs in the orbit $\Omega_2$. We also use the notation $U_1(\mO'), U_2(\mO')$, if $\mO' = \gamma_i(\mO)$. Remark that $U_1$ and $U_2$ depend on the embedding $\gamma_{i}$, and that $U_i(\mO)$ only depends on the choice of the orbit $\Omega_i$.

Now consider any isomorphism 
\begin{equation}
\beta: \mQ\ \mapsto\ \mQ''
\end{equation}
with $\mQ''$ a $\mQ(4,q)$-quadrangle in $\Gamma$, such that $\beta(\mO)$ is a subtended ovoid in $\mQ''$. Then $\beta^{-1}([\infty]) \in U_1(\mO) \cup U_2(\mO)$. \\

\begin{remark}
{\rm
Later, in section \ref{later}, we will see that $\vert U_1(\mO) \vert \equiv 0\mod{2}$ if $\sigma^2 = \id$. }
\end{remark}

\subsection{Intrinsic way to recover the embedding}
\label{intr}

Given a fixed element of  $\texttt{Emb}(\mQ(4,q),\Gamma)$, there is an intrinsic way to recover the set $\mO(\mQ)$, once we know at least one element of $\mO(\mQ)$. \\

Put $C := \Aut(\Gamma)_{\mQ}$. 
Each line in $[\infty]^{\perp} \cap \mQ$ is an axis of symmetry (see section \ref{recapor}); denote the group generated by the symmetries defined by these lines by $L = L_{[\infty]}$. Note that 
$L \leq C$.

Now suppose $x \I [\infty]$ is an arbitrary point, and let $y, z$ be external points to $\mQ$ which are contained in $x^{\perp}$. If they are collinear, obviously 
there is an element of $L$ mapping $y$ to $z$. Suppose they are not collinear, so that $yx =: U$ and $zx =: V$ are different lines. Now consider a $\mQ(4,q)$-subGQ $\widetilde{\mQ}$
containing $[\infty], U$ and $V$, such that $\mQ \cap \widetilde{\mQ}$ is a grid $\mG$ with parameters $(q,1)$ (note that each of the $q - 1$ $\mQ(4,q)$-subGQs distinct from $\mQ$
containing $[\infty], U$ and $V$ meets $\mQ$ in such a grid (since they do not share all the lines on $x$)). Let $M$ and $N$ be different lines in $\mG \cap [\infty]^\perp$. 
Let $L(M,N)$ denote the subgroup of $L$ generated by the symmetries about $M$ and $N$; this group is isomorphic to $\mathbf{SL}_2(q)$ by \cite{two,SFGQ}, stabilizes $\mQ$ and $\widetilde{\mQ}$, 
and  acts transitively on the points of $\mQ$, respectively $\widetilde{\mQ}$, which are not contained in $\mG$. In particular, there is an element in $L(M,N)$ 
which sends $y$ to $z$. 

Since $L(M,N) \leq L$, and since
$L$ acts transitively on the points of $[\infty]$, it now easily follows that $L$ acts transitively on the points of $\Gamma$ outside $\mQ$. Whence it acts 
transitively on the subtended Kantor-Knuth ovoids in $\mQ$. So if $\mO_x$ is such an ovoid, we have that
\begin{equation}
{\mO_x}^{C} = {\mO_x}^L.
\end{equation}

Remark that the result  holds for any element $\mQ$ in $\Omega_ 1 \cup \Omega_2$. \\

By the preceding subsection, we have the following result. 

\begin{theorem}[Intrinsic Kantor-Knuth embedding]
Consider an embedding 
\begin{equation}
\gamma:\ \mQ\ \mapsto\ \mQ',
\end{equation}
where $\mQ \cong \mQ(4,q)$ and $\mQ'$ is a subGQ of a nonclassical Kantor-Knuth quadrangle $\Gamma$ of order $(q,q^2)$. Let $\mO$ be any Kantor-Knuth ovoid in $\mQ$ 
such that $\gamma(\mO)$ is a subtended ovoid in $\mQ'$. Let $[\infty]$ be the special line in $\Gamma$, and put $\gamma^{-1}([\infty]) = U \in U_1(\mO) \cup U_2(\mO)$. 
Then 
\begin{equation}
\gamma^{-1}\Big(\mO(\mQ')\Big)\ =\ \mO^{L_U}.
\end{equation}
\eop
\end{theorem}

We still need to obtain the result mentioned at the end of the previous section, that is, that $U_1^{\Aut(\mQ)_\mO} \ne U_2^{\Aut(\mQ)_\mO}$. We use the notation of that section. By the results of the present section, 
it follows that $\vert \mO(\mQ'_i) \vert = \vert \mO^{L_{U_i}} \vert$; but as $i$ differs, this quantity also changes by a factor $2$. It follows that there is no 
element in $\Aut(\mQ)_\mO$ which maps $U_1$ to $U_2$.



\medskip
\subsection{Remark: configuration of special lines}

We keep using the notation of the previous section.

Relative to one given embedding $\gamma$, one might wonder how the configuration of special lines corresponding to the embedded Kantor-Knuth ovoids looks like. The answer is very simple (and we take $U_1 := \gamma^{-1}([\infty])$).\\

Consider the set $\bigcup_{\mO \in \gamma^{-1}(\mO(\mQ'))}\Big( U_1(\mO) \cup U_2(\mO) \Big)$ with $\mQ' \in \Omega_1$; then by the fact that $L_{U_1}$ acts on $\gamma^{-1}(\mO(\mQ'))$, we have:
\begin{equation}
\bigcup_{\mO \in \gamma^{-1}(\mO(\mQ'))}\Big( U_1(\mO) \cup U_2(\mO) \Big) \ =\ U_1^{\perp}.
\end{equation}

Also, $\bigcup_{\mO \in \gamma^{-1}(\mO(\mQ'))}  U_2(\mO)  =  U_1^{\perp} \setminus \{ U_1 \}$. 


\medskip
\subsection{An equivalence relation on the Kantor-Knuth orbit}

Let $\F_q$ be a finite field for $q$ an odd prime power, and let $\mQ(4,q)$ be an orthogonal quadrangle in projective $4$-space.  
Let $\Omega(q,\sigma)$ be the complete $\Aut(\mQ(4,q))$-orbit of Kantor-Knuth ovoids of type $\sigma \in \Aut(\F_q)^\times$. Take any element $\mO \in \Omega(q,\sigma)$, pick $U_1 \in U_1(\mO)$, $U_2 \in U_2(\mO)$, and define 
\begin{equation}
\begin{cases}
\mA_1\ :=\ \{ (\eta(\mO),\eta(U_1))\ \vert\ \eta \in \Aut(\mQ(4,q)) \} \\
\mA_2\ :=\ \{ (\eta(\mO),\eta(U_2))\ \vert\ \eta \in \Aut(\mQ(4,q)) \}.
\end{cases}
\end{equation}

By the above, $\mA_1$ and $\mA_2$ are different (so disjoint) subsets of $\Omega(q,\sigma) \times \mL$, where $\mL$ is the line set of $\mQ(4,q)$. Note that the natural projection
\begin{equation}
\mA_1 \cup \mA_2\ \mapsto \mL
\end{equation}
is surjective.

Obviously, $\mA_1$ corresponds to subGQs in $\Omega_1$ in the Kantor-Knuth GQ $\Gamma(q,\sigma)$, and $\mA_2$ corresponds to subGQs in the other orbit. 

Now define an equivalence relation $\sim_\epsilon$ on $\mA_i$, $i = 1,2$, as follows: $(\mO,U) \sim_\epsilon (\mO',U')$ if and only if $U = U'$ and $\mO \in {\mO'}^{L_{U'}}$. Note that this indeed defines an equivalence relation on both $\mA_1$ and $\mA_2$. 

\begin{observation}
Each equivalence class of $\sim_{\epsilon}\ \subset \mA_i$ exactly singles out the subtended ovoids in an isomorphic image of $\mQ(4,q)$ in $\Gamma = \Gamma(q,\sigma)$, and 
$i$ determines in which orbit the isomorphic image of $\mQ(4,q)$ is situated. \eop 
\end{observation}






\medskip
\subsection{Embedding theorem}

As of now, we will denote objects of $\texttt{Emb}(\mQ(4,q),\Gamma)$ sometimes as $\Big[ \mQ(4,q),\gamma,\mQ \Big]$.


\begin{theorem}[Embedding theorem]
\label{embthm}
Consider two objects $A := \Big[ \mQ(4,q),\gamma,\mQ \Big]$ and $B := \Big[ \mQ(4,q),\gamma',\mQ' \Big]$ in $\texttt{Emb}(\mQ(4,q),\Gamma)$. Then there exists 
an isomorphism $\kappa = (\alpha,\beta): A \mapsto B$ if and only if $\mQ$ and $\mQ'$ belong to the same $\Aut(\Gamma)$-orbit.  
In particular, if $A = B$, $\Aut(A) \cong \Aut(\mQ)_{\mO(\mQ)}$. 
\end{theorem}

{\em Proof}.\quad 
Let $\kappa = (\alpha,\beta): A \mapsto B$ be an isomorphism; then $\beta(\mO(\mQ)) = \mO(\mQ')$, so that if $U$ is any special line for $\mO \in \mO(\mQ)$, then 
$\vert \mO(\mQ) \vert = \vert \mO^{L_U} \vert = \vert \mO^{L_{\beta(U)}} \vert = \vert \mO(\mQ') \vert$. This can only happen if $\mQ$ and $\mQ'$ are in the same 
$\Gamma$-orbit. 

The other direction is obvious.\\

Now let $A = B$, and consider $\Aut(A)$. We are looking for pairs of automorphisms $(\alpha,\beta)$ such that the diagram below commutes:

\begin{center}
\begin{tikzpicture}[>=angle 90,scale=3.5,text height=2.0ex, text depth=0.45ex]
\node (a0) at (0,3) {$\mQ$};
\node (a0') [right=of a0] {};
\node (a1) [right=of a0'] {$\mQ$};

\node (b0) [above=of a0] {$\mQ(4,q)$};
\node (b1) [above=of a1] {$\mQ(4,q)$};

\draw[->,font=\scriptsize,thick,>=angle 45,orange]

(a0) edge node[above] {$\beta$} (a1)
(b0) edge node[above] {$\alpha$} (b1);



\draw[->,font=\scriptsize,thick]

(b1) edge node[right] {$\gamma$} (a1);

\draw[->,font=\scriptsize,thick]
(b0) edge node[left] {$\gamma$} (a0);

\end{tikzpicture}
\end{center}
with $\beta(\mO(\mQ)) = \mO(\mQ)$. 

Obviously, such pairs are completely and uniquely determined by the automorphisms $\beta$ of $\mQ$ which stabilize the set $\mO(\mQ)$. This is what we had to prove.
\eop \\


\medskip
\subsection{Isomorphism classes of $\texttt{Emb}(\mQ(4,q),\Gamma)$}

Note that by Theorem \ref{embthm}, $\texttt{Emb}(\mQ(4,q),\Gamma)$ only has two isomorphism classes of objects. For essentially the same reason, 
the same is true for $\widetilde{\texttt{Emb}}(\mQ(4,q),\Gamma)$.



\medskip
\subsection{The categories for different Kantor-Knuth quadrangles}

Now let $\Gamma$ and $\Gamma'$ be nonisomorphic Kantor-Knuth quadrangles of the same order $(q,q^2)$ | say with respective (different) associated $\F_q$-automorphisms 
$\sigma$ and $\sigma'$. Question: are $\texttt{Emb}(\mQ(4,q),\Gamma)$ and $\texttt{Emb}(\mQ(4,q),\Gamma')$ isomorphic?  

Consider an object $A$ in $\texttt{Emb}(\mQ(4,q),\Gamma)$; then by Theorem \ref{embthm}, $\Hom(A,A) = \Aut(A)$ can be identified with the automorphism 
group $\Aut(\mQ)_{\mO(\mQ)}$ of $\mQ$ which fixes $\mO(\mQ)$ (this group automatically fixes $[\infty]$; for if some element $\varphi$ would move $[\infty]$, it is straightforward to see that $L_{[\infty]}$ and $L_{[\infty]^\varphi}$ generate a group which acts transitively on the set of all lines, and this is sufficient information to force a contradiction: {the set of special points of the ovoids in $\mO(\mQ)$ coincides with the set of points incident with $[\infty]$, so $[\infty]$ must be fixed by $\Aut(\mQ)_{\mO(\mQ)}$}).  
As we have seen, this group acts transitively on the elements of $\mO(\mQ)$. Note also that $\hom(A,A)$ is independent of the choice of $\mQ$ in $\Omega_1$ or $\Omega_2$. 
Since $\Gamma$ and $\Gamma'$ are not isomorphic, $\sigma$ and $\sigma'$ have a different order.
Now if one of $\sigma$, $\sigma'$ would be an involution, and the other is not, then 
$\Hom(A,A)$ and $\Hom(B,B)$ with $B$ in $\texttt{Emb}(\mQ(4,q),\Gamma')$ are not isomorphic as abstract groups, since their sizes are different (as for instance 
$\vert \Hom(A,A) \vert = \vert \mO(\mQ) \vert \cdot \vert \Aut(\mO) \vert$, with $\mO \in \mO(\mQ)$).
So the categories cannot be isomorphic. If we suppose that both $\sigma$ and $\sigma'$ are involutions, then  
$\Hom(A,A)$ and $\Hom(B,B)$ do have the same size, so could be isomorphic in principle.\\

The following easy theorem reveals one crucial difference between the action of $\Aut(\mQ)_{\mO(\mQ)}$ on $\mO(\mQ)$ in case $\Gamma$ varies. 
First consider any nonclassical Kantor-Knuth ovoid $\mO$ in $\mQ(4,q)$ with special automorphism $\sigma$, and let $u$ be its special point; furthermore, let $v \in \mO \setminus \{u\}$. 
Then by definition, $H(u,v)$ is the subgroup of $\Aut(\mQ(4,q))_\mO$ which fixes $u$ and $v$ linewise; we call this group {\em kernel group} or {\em kernel} of $\mO$, and it is 
isomorphic to $\F_{\sigma}^\times$, where $\F_{\sigma}$ is the subfield of $\F_q$ fixed elementwise by $\sigma$ (see e. g. \cite{bloem} for details). 

\begin{theorem}
\label{kern}
Let $\Gamma$ be a nonclassical Kantor-Knuth GQ of order $(q,q^2)$, $\mQ \in \Omega_1 \cup \Omega_2$ be any $\mQ(4,q)$-subGQ, and $\mO \in \mO(\mQ)$. 
Then $\Aut(\Gamma)_{\mQ,\mO(\mQ)} = \Aut(\Gamma)_\mQ$ generates the kernel of $\mO$. By ``generates the kernel,''  we mean that the kernel group $H(u,v)$, with $u = [\infty] \cap \mO$ and $v \in \mO \setminus \{ u\}$,  of any ovoid $\mO \in \mO(\mQ)$ is 
contained in $\Aut(\Gamma)_{\mQ,\mO(\mQ)}$. (That is: $G(u,v) \cap \Aut(\Gamma)_{\mQ,\mO(\mQ)} \geq H(u,v)$, with $G(u,v)$ the kernel of the TGQ $\Gamma$.)
\end{theorem}

{\em Proof}.\quad
Let $\mO$ be as in the statement, let $u = [\infty] \cap \mO$ and let $v \in \mO \setminus \{ u\}$ be arbitrary. Consider the kernel  $G(u,v)$ of the TGQ $\Gamma$ with respect to $u$ and $v$ (cf. \cite[\S 3.4]{TGQ}). We have that $G(u,v)$ is isomorphic to $\F_{\sigma}^\times$ (cf. \cite[4.7.3]{TGQ}).  
Then each element $\nu \in G(u,v)$ fixes each line of $\mQ$ incident with $u$ or $v$, so applying  \cite[Lemma 4.2.5]{SFGQ} on $\mQ \cap \mQ^\nu$, we 
conclude that $\mQ$ is fixed by $\nu$. So $G(u,v)$ faithfully induces an automorphism group of $\mQ$ (faithfully because a nontrivial automorphism of $\Gamma$ that fixes 
$\mQ$ elementwise, cannot fix lines outside $\mQ$). Now let $w$ be any point in $\{ u, v \}^{\perp} \setminus \mQ$; then $G(u,v)$ fixes $\mO_w = w^{\perp} \cap \mQ$. 
Suppose w. l. o. g. that $\mO = \mO_x$, with $x$ a point of $\Gamma \setminus \mQ$ subtending $\mO_x$. Then $x \in \{ u, v\}^{\perp}$, so if we put $w = x$, 
we deduce that $G(u,v)$ fixes $\mO_x$. Whence $G(u,v)$ induces the entire kernel group $H(u,v) \cong G(u,v)$ of $\mO_x$. The statement of the theorem follows. 
\eop \\

\begin{corollary}
We have that $\Aut(\mQ)_{\mO(\mQ)}$ generates the kernel of any such $\mO$. 
\end{corollary}

{\em Proof}.\quad
Immediate as $\Aut(\Gamma)_{\mQ,\mO(\mQ)} \Big/ N \leq \Aut(\mQ)_{\mO(\mQ)}$, with $N$ the kernel of the action of $\Aut(\Gamma)_{\mQ,\mO(\mQ)}$ on $\mQ$. 
\eop \\

Even if $\Hom(A,A)$ and $\Hom(B,B)$ would be isomorphic, the corresponding actions are not, as the following theorem shows.

\begin{theorem} 
Let $\Gamma$ and $\Gamma'$ be nonisomorphic Kantor-Knuth GQs of order $(q,q^2)$, and let $\mQ \leq \Gamma$ and $\mQ' \leq \Gamma'$ be $\mQ(4,q)$-subGQs. Then the actions $\Aut(\mQ)_{\mO(\mQ)} \curvearrowright \mO(\mQ)$ and $\Aut(\mQ)_{\mO(\mQ')} \curvearrowright \mO(\mQ')$ are not equivalent.
\end{theorem}

{\em Proof}.\quad
Since $\Gamma$ and $\Gamma'$ are defined by different automorphisms of $\Aut(\F_q)$, the corresponding kernels are different subfields of $\F_q$. So if $\mO \in \mO(\mQ)$ 
and $\mO' \in \mO(\mQ')$, they define different isomorphism classes of kernel groups. This implies the theorem. \eop \\

\medskip
\subsection{A weird observation}

Let $U$ be any line of $\mS = \mQ(4,q)$, and define $L = L_U$ as the group generated by all symmetries with axis contained in $U^{\perp}$. Let $H(L)$ be the kernel group,  with respect to the points $u$ and $v$,  with $u \I U$ and $v \not\sim u$,  generated 
by $L$ in the sense of the previous sections: it is $G(u,v) \cap L$ with $G(u,v)$ the kernel of the TGQ $\mS$ with translation point $(\infty) = u$ (and with respect to $v$).  

We will prove the following  result without any calculation, only using embedding theory. 

\begin{theorem}
If $q$ is odd, $H(L)$ is a subgroup of $\F_p^{\times}$, where $q$ is a power of the prime $p$. 
\end{theorem}

{\em Proof}.\quad
Suppose $q$ is odd (and fixed), and consider any embedding $\gamma: \mQ(4,q) \mapsto \mQ$, with $\mQ$ a subquadrangle of a nonclassical Kantor-Knuth GQ $\Gamma$, and 
$\gamma(U)$ the line $[\infty]$ of $\Gamma$.  Since $L_U = \gamma^{-1}(L_{[\infty]})$, and as $L_U$ and 
$H(L_U)$ are independent of $\Gamma$, we can consider all such embeddings. For a fixed $\gamma$, we have seen that the kernel group generated by $L_{[\infty]}$ is contained 
in the kernel of the Kantor-Knuth TGQ (once we have chosen a suitable translation point).  

It easily follows that $H(L_U)$ is contained in the intersection of the kernel groups obtained by letting $\gamma$ (and so also $\Gamma$, and the automorphism $\sigma \in \Aut(\F_q)$) vary. So $H(L_U)$ must be a subgroup of $\F_p^{\times} = \bigcap_{\sigma \in \Aut(\F_q)}\F_\sigma^{\times}$. 
\eop \\

\medskip
\section{The automorphism group via the good orbit?}
\label{later}

{
Let $\mQ \cong \mQ(4,q)$ be any subGQ of order $(q,q)$ in a given GQ $\Gamma$ isomorphic to $\mQ(5,q)$; then $\mE$, the geometry of subtended ovoids and rosettes of ovoids, is $2$-covered by $\mA = \Gamma \setminus \mQ$. Now it is a simple fact that 
\begin{equation}
\Aut(\Gamma)_{\mQ}/N\ \cong\ \Aut(\mE), 
\end{equation}
with $N$ the kernel of the action of $\Aut(\Gamma)_{\mQ}$ on $\mQ$, and that if $x$ is a point in $\Gamma \setminus \mQ$, then $\Aut(\mQ)_{\mO_x}$  is induced by $\Aut(\Gamma)_{\mQ,x}$. 

So we ask ourselves whether the same is true if we replace $\Gamma$ by a nonclassical Kantor-Knuth quadrangle, and $\mQ$ taken in the appropriate subGQ-orbit. 
Since the size of the subGQ-orbit $\Omega_1$ is rather small compared to the size of $\Omega_2$, and since each subGQ $\mQ \in \Omega_1$ is $2$-covered by $\Gamma \setminus \mQ$, it is a natural question to ask how big a part 
of the full automorphism group (in $\mQ$) of an element $\mO \in \mO(\mQ)$ is recovered by its automorphism group induced by $\Aut(\Gamma)_{\mQ}$. Since the special line $[\infty]$ is fixed by $\Aut(\Gamma)$, we will be at least a factor $2$ off as it is not fixed by $\Aut(\mQ)_\mO$. Still, as we will see below in a simple calculation, the automorphism group induced by $\Aut(\Gamma)_\mQ$ will be very large in comparison to $\Aut(\mO)$: indeed, it will turn out to have half the size of $\Aut(\mQ)_\mO$.   }\\

Keep using the notation as in the previous section. Since $\vert \Omega_1 \vert = 2q^2$, it is easy to observe that 
the size of $\Aut(\Gamma)_{\mO_x} = {\Big(\Aut(\Gamma)_{\mQ}\Big)}_{\mO_x}$ is
\begin{equation}
\label{bound1}
\frac{(q + 1)(q - 1)^2q^6\delta h}{2q^2\cdot \Big(q^2(q + 1)(q - 1)/2\Big)} = (q - 1)q^2\delta h. 
\end{equation}

In $\mQ$, this group induces a group of size $(q - 1)q^2\frac{\delta}{2}h$. 

By a calculation in \cite{PW}, the size of the complete automorphism group of $\mO_x$ in $\Aut(\mQ)$ is $(q - 1)q^2\cdot \delta h$. It is known 
that $\Aut(\mO_x) := \Aut(\mQ)_{\mO_x}$ fixes some point $u \in \mO_x$ (recall section \ref{tranov}); let $U$ be any line incident with $u$, and put $\ell := \vert U^{\Aut(\mQ)_{\mO_x}} \vert$.

As $\Aut(\mQ)_U$
is isomorphic to $\mathbf{P\Gamma L}_2(q)$ in its action on $U$, we know that $\vert \Aut(\mQ)_U/M \vert = (q + 1)q(q - 1)h$, with $M$
the kernel of the action of $\Aut(\mQ)_U$ on $U$ (the latter seen as a point set). So
\begin{equation}
\Big \vert \Aut(\mO_x)\Big \vert \Big/ \Big(\Big\vert U^{\Aut(\mO_x)} \Big\vert \cdot \Big\vert M \cap \Aut(\mO_x)\Big\vert\Big) \ \ \mbox{divides}\  \ q(q - 1)h.
\end{equation}

{
Here, 
\begin{equation}
\frac{\Big\vert \Aut(\mO_x) \Big\vert}{\Big\vert U^{\Aut(\mO_x)} \Big\vert}\ =\ \Big\vert \Aut(\mO_x)_U \Big\vert,  
\end{equation}
and we have to mod out the trivial action of $M \cap \Aut(\mO_x)$.} Note that ${\Aut(\mO_x)}_U$ is contained in $\Aut(\mQ)_{U,v}$, with $v = \mO \cap U$. 

We will use the following fact below:
\begin{quote}
Fact: let $Y$ be a full grid in a $\mQ(4,q)$ with $q$ odd. Then the elementwise stabilizer of $Y$ in $\Aut(\mQ(4,q))$ has size $2$. 
\end{quote}

Consider a Sylow $2$-subgroup $S_2$ of $M \cap \Aut(\mO_x)$. 
Then $S_2$ fixes $U$ pointwise. Let $v \I U$ be any point incident with $U$; then by the fact that $q$ is odd, 
$S_2$ has to fix at least one other line $V$ incident with $v$ besides $U$. Each such line meets $\mO_x$ in a point, which is necessarily fixed by $S_2$. As 
no two points of $\mO_x$ are collinear, it easily follows that the fixed elements structure of the group $S_2$ is a full grid 
of order $(q,1)$ if $S_2 \ne \{ \id \}$. 

By the aforementioned fact, It follows that $\delta/2$ divides  $\ell = \vert U^{\Aut(\mQ)_{\mO_x}} \vert$.

If $q$ is a square, and $\sigma^2 = \id \ne \sigma$, then $\delta = 4$, so each such orbit has size $\ell$ divisible by $2$. In case $U = [\infty]$, $\ell = 2$  is also the maximal size of such an orbit, since 
we know that the full automorphism group of $\mO_x$ inside $\mQ$ has size  $(q - 1)q^2\delta h$, and since the action of $\Aut(\mQ)_{\mO_x}$ on $\mO_x$ is 
faithful. {If $\sigma^2 \ne \id$, then $\delta = 2$, so the divisibility condition becomes vacuous. 
In any case, when $U = [\infty]$ we have that $\ell \leq 2$, since ${\Big(\Aut(\Gamma)_\mQ\Big)}_{\mO_x}$ induces an automorphism group of $\mO_x$ inside $\mQ$ which fixes $[\infty]$.} 


\subsection*{Comparing the different sizes}

As we have seen, for $\mO_x$ as above, we obtain that its stabilizer inside $\Big(\Aut(\Gamma)\Big)_{\mQ}$ has size $(q - 1)q^2\frac{\delta}{2}h$, just a factor $2$ 
shy of the size of its full automorphism group.




\end{document}